\documentclass[a4paper, 11pt]{amsart}

\def\siecle#1{\textsc{\romannumeral #1}\textsuperscript{e}~si\`{e}cle}

\usepackage{csquotes}

\usepackage[margin=3cm]{geometry}

\usepackage{amsmath,amssymb,amsfonts,amsthm}
\usepackage[all]{xy}
\usepackage[latin1]{inputenc}
\usepackage{enumerate}
\usepackage[colorlinks=true,linkcolor=red,citecolor=blue]{hyperref}
\usepackage{tikz}
\usepackage{fourier}

\usepackage[frenchb]{babel}

\theoremstyle{theorem}
\newtheorem*{thm}{Th\'{e}or\`{e}me}
\newtheorem*{conj}{Conjecture}

\theoremstyle{remark}
\newtheorem*{exa}{Exemple A}
\newtheorem*{exb}{Exemple B}

\theoremstyle{definition}
\newtheorem*{defi}{D\'{e}finition}

\renewcommand{\leq}{\leqslant}
\renewcommand{\geq}{\geqslant}
\newcommand{\ZZ}{\mathbb{Z}}
\newcommand{\QQ}{\mathbb{Q}}
\newcommand{\RR}{\mathbb{R}}
\newcommand{\CC}{\mathbb{C}}

\renewcommand{\to}{\rightarrow}
\newcommand{\To}{\longrightarrow}
\renewcommand{\H}{\operatorname{H}}
\newcommand{\M}{\operatorname{M}}
\newcommand{\B}{\mathrm{B}}
\newcommand{\dR}{\mathrm{dR}}

\usepackage{microtype}

\begin{document}

\title{Motifs : un tour d'horizon}

%
\author{Cl\'{e}ment Dupont}

\address{Institut Montpelli\'erain Alexander Grothendieck, Universit\'e de Montpellier, CNRS,
Montpellier, France}
\email{clement.dupont@umontpellier.fr}

%

\maketitle

Ce texte est paru dans la \emph{Gazette de la Soci\'{e}t\'{e} Math\'{e}matique de France} en octobre 2023. Il constitue une introduction (partielle et partiale) \`{a} la th\'{e}orie des motifs, un programme initi\'{e} par Grothendieck et qui structure aujourd'hui la g\'{e}om\'{e}trie et l'arithm\'{e}tique des vari\'{e}t\'{e}s alg\'{e}briques.

\section{Aux sources des motifs}

	Nous d\'{e}crivons ici une partie de la g\'{e}n\'{e}alogie, multiple et diverse, de la notion de motif, en mettant l'accent sur le r\^{o}le important des conjectures de Weil. 
	Ces conjectures concernent les fonctions z\^{e}ta des vari\'{e}t\'{e}s alg\'{e}briques d\'{e}finies sur un corps fini et trouvent elles-m\^{e}mes leur source dans les travaux de Riemann sur la r\'{e}partition des nombres premiers. 
	La th\'{e}orie des motifs telle que pens\'{e}e par Grothendieck est un prolongement de travaux classiques sur la cohomologie des vari\'{e}t\'{e}s alg\'{e}briques et la g\'{e}om\'{e}trie \'{e}num\'{e}rative, que nous introduisons.

	\subsection{Fonction z\^{e}ta, fonctions $L$}\label{par: zeta}

		Riemann introduit en 1859 la \emph{fonction z\^{e}ta}
		$$\zeta(s) = \sum_{n=1}^{+\infty}\frac{1}{n^s} \quad (s\in\mathbb{C}, \operatorname{Re}(s)>1),$$
		qu'il prolonge en une fonction m\'{e}romorphe sur tout le plan complexe, et pour laquelle il d\'{e}montre une \'{e}quation fonctionnelle reliant $\zeta(s)$ et $\zeta(1-s)$. Euler avait d\'{e}j\`{a} calcul\'{e} en 1735 :
		\[
		\zeta(2)=\frac{\pi^2}{6} \; , \; \zeta(4)=\frac{\pi^4}{90} \; , \; \zeta(6)=\frac{\pi^6}{945} \;,\;\ldots,
		\]
		et en g\'{e}n\'{e}ral $\zeta(2k)\in\QQ\pi^{2k}$ pour tout $k\in\mathbb{N}^*$.
		
		La factorisation en produit de nombres premiers donne une \'{e}criture de la fonction z\^{e}ta, d\'{e}j\`{a} remarqu\'{e}e par Euler et appel\'{e}e \emph{produit eul\'{e}rien} :
		$$\zeta(s) = \prod_{p \mbox{ \footnotesize{premier}}} \frac{1}{1-p^{-s}} \ \cdot$$
		Riemann en d\'{e}duit que les z\'{e}ros complexes de la fonction z\^{e}ta contr\^{o}lent la r\'{e}partition des nombres premiers, et formule au passage une conjecture toujours ouverte, l'\emph{hypoth\`{e}se de Riemann} : ces z\'{e}ros (hormis les z\'{e}ros \enquote{triviaux} $s=-2,-4,-6,\ldots$) devraient \^{e}tre situ\'{e}s sur la droite $\operatorname{Re}(s)=\frac{1}{2}$.
				
		Vingt ans avant Riemann, son mentor Dirichlet avait d\'{e}j\`{a} us\'{e} des m\'{e}thodes de l'analyse (r\'{e}elle) pour prouver le \emph{th\'{e}or\`{e}me de la progression arithm\'{e}tique} : pour des entiers $a,d\geq 1$ premiers entre eux, il existe une infinit\'{e} de nombres premiers congrus \`{a} $a$ modulo $d$. Sa preuve repose sur les \emph{fonctions $L$},
		$$L_\chi(s) = \sum_{n=1}^{+\infty}\frac{\chi(n)}{n^s},$$
		associ\'{e}es \`{a} des \emph{caract\`{e}res de Dirichlet} $\chi:\mathbb{Z}\to \mathbb{C}^*$. 
		
		
		\`{A} partir de la fin du \siecle{19} apparaissent de nombreuses \enquote{fonctions $L$} (ou fonctions $\zeta$) dans l'esprit de celles de Dirichlet, associ\'{e}es \`{a} divers objets de nature arithm\'{e}tique : corps de nombres (Dedekind), repr\'{e}sentations galoisiennes (E. Artin), formes modulaires (Hecke), etc. Ces fonctions ont des caract\'{e}ristiques (souvent conjecturales !) communes : produit eul\'{e}rien, prolongement analytique, \'{e}quation fonctionnelle, valeurs int\'{e}ressantes aux entiers, hypoth\`{e}se de Riemann.
		
		Une des sources de la th\'{e}orie des motifs 
		est la recherche d'un cadre pour structurer les propri\'{e}t\'{e}s des fonctions $L$ issues de l'arithm\'{e}tique, qui sont toutes des cas particuliers de fonctions $L$ associ\'{e}es aux motifs.

	\subsection{Co\"{i}ncidences ?}
	
		Le terme de \emph{motif}, invent\'{e} par Grothendieck, \'{e}voque une phrase musicale qu'on	 retrouve, sous diff\'{e}rentes incarnations, \`{a} plusieurs endroits d'une \oe uvre. En peinture, le motif est le sujet d'une toile, qui en est une repr\'{e}sentation du point de vue d'un artiste. Un motif au sens math\'{e}matique est un objet associ\'{e} \`{a} une vari\'{e}t\'{e} alg\'{e}brique et qui s'incarne en des invariants de natures diff\'{e}rentes, expliquant certaines \enquote{co\"{i}ncidences} entre ces invariants. C'est donc aussi le motif, au sens de raison d'\^{e}tre, de ces incarnations diff\'{e}rentes d'un ph\'{e}nom\`{e}ne commun.
		
		Une de ces co\"{i}ncidences a trait \`{a} la notion de genre. Pour une surface de Riemann 
		compacte $X$, on a l'\'{e}galit\'{e} entre les deux quantit\'{e}s suivantes, appel\'{e}es \emph{genre} de $X$ (voir \cite{fle, popescupampugenus} pour des analyses historiques de la notion de genre).
		\begin{enumerate}
		\item[(a)] Un invariant \emph{topologique} : le \enquote{nombre de trous} de $X$, c'est-\`{a}-dire le nombre de tores (\enquote{donuts}) qu'on doit recoller entre eux pour former $X$. C'est-aussi le nombre maximal de coupures (donn\'{e}es par des courbes continues, ferm\'{e}es simples deux \`{a} deux disjointes) qu'on peut effectuer sur $X$ sans la d\'{e}connecter. La figure suivante\footnote{Les figures de cet article ont \'{e}t\'{e} r\'{e}alis\'{e}es par Anthony Genevois.} 
		montre une surface de Riemann compacte de genre $2$ avec un syst\`{e}me maximal de coupures, et une repr\'{e}sentation de la surface d\'{e}coup\'{e}e.
		\end{enumerate}
		\begin{figure}[h!]
		\includegraphics[scale=0.29]{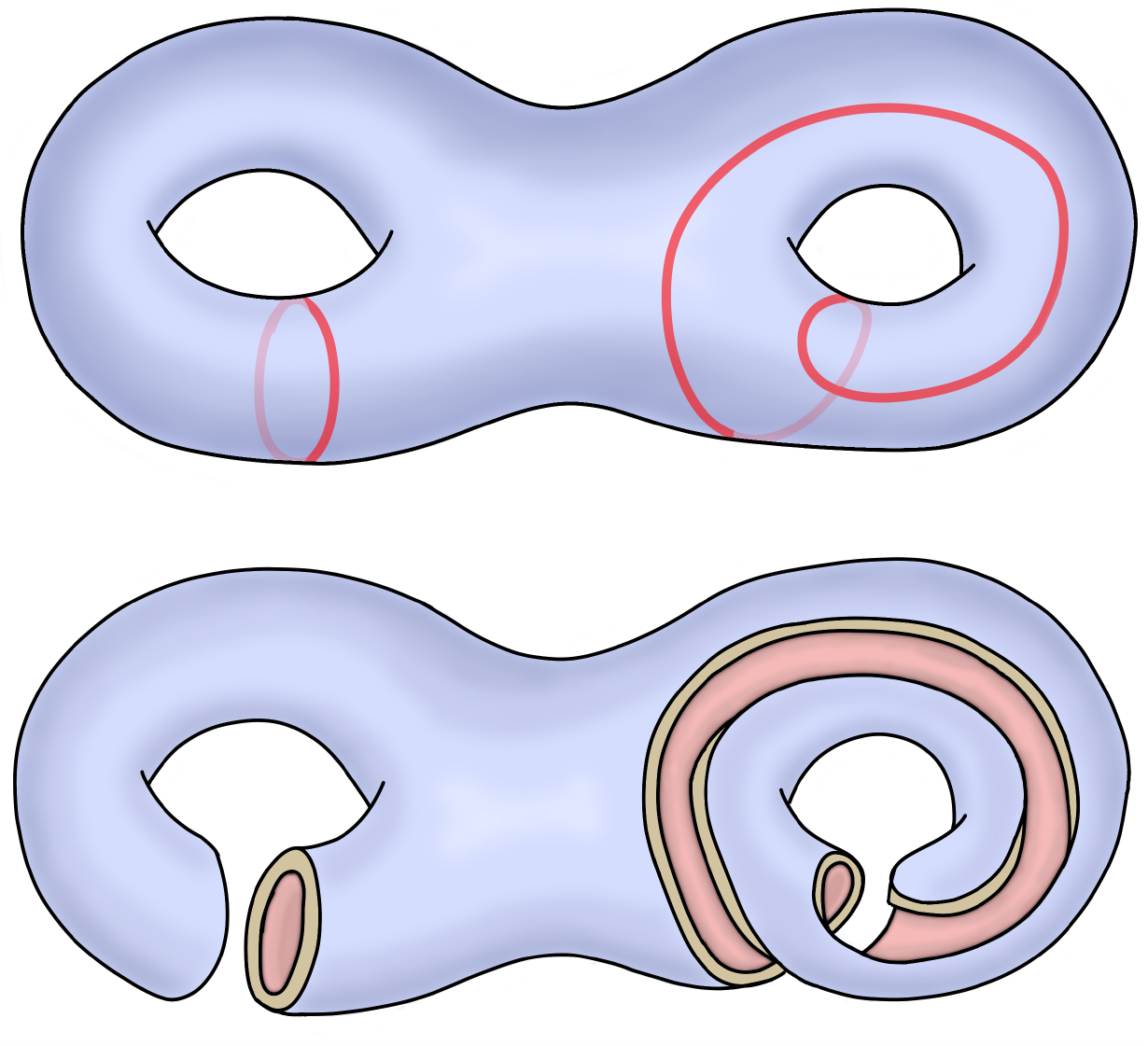}
		\end{figure}
		\begin{enumerate}
		\item[(b)] Un invariant \emph{analytique} : la dimension de l'espace vectoriel complexe des $1$-formes diff\'{e}rentielles holomorphes sur $X$ (donn\'{e}es dans une coordonn\'{e}e locale $z$ par une formule $f(z)\operatorname{d}\!z$ avec $f$ holomorphe), classiquement appel\'{e}es \enquote{diff\'{e}rentielles de premi\`{e}re esp\`{e}ce}. 
		\end{enumerate}

		On verra au \S\ref{par: conjectures de weil} une troisi\`{e}me incarnation du genre qui a, via les conjectures de Weil, influenc\'{e} les d\'{e}veloppements de la th\'{e}orie des motifs.
		
		\begin{enumerate}
		\item[(c)] Un invariant \emph{arithm\'{e}tique}, li\'{e} au nombre de solutions d'un syst\`{e}me d'\'{e}quations polynomiales correspondant \`{a} une courbe alg\'{e}brique d\'{e}finie sur un corps fini.
		\end{enumerate}
		
		Les  \enquote{co\"{i}ncidences} dont la notion de motif est cens\'{e}e rendre compte se formalisent via la cohomologie des vari\'{e}t\'{e}s alg\'{e}briques.
		
	\subsubsection{Vari\'{e}t\'{e}s alg\'{e}briques}
		
		On s'int\'{e}resse aux vari\'{e}t\'{e}s alg\'{e}briques d\'{e}finies sur un corps $k$, comme les vari\'{e}t\'{e}s affines (resp. projectives), d\'{e}finies dans l'espace affine $\mathbb{A}^n$ (resp. projectif $\mathbb{P}^n$) par un syst\`{e}me d'\'{e}quations donn\'{e}es par des polyn\^{o}mes en $n$ variables (resp. des polyn\^{o}mes homog\`{e}nes en $n+1$ variables) \`{a} coefficients dans $k$. Pour une vari\'{e}t\'{e} alg\'{e}brique $X$ et un corps $K$ contenant $k$, on a l'ensemble $X(K)$ de ses $K$-points rationnels, qui est l'ensemble des solutions dans $K$ du syst\`{e}me d'\'{e}quations.
		
		
		\begin{exa}
		Consid\'{e}rons $X=\mathbb{A}^1\setminus \{0\}$, la droite affine priv\'{e}e du point $0$, d\'{e}finie sur $k=\QQ$. C'est une vari\'{e}t\'{e} alg\'{e}brique affine qui peut \^{e}tre d\'{e}crite par l'\'{e}quation $xy=1$ dans le plan affine $\mathbb{A}^2$.  On a $X(\CC)=\CC^*$.
		\end{exa}
		
		
		
		\begin{exb}
		Consid\'{e}rons la vari\'{e}t\'{e} alg\'{e}brique d\'{e}finie sur $k=\QQ$ par une \'{e}quation de la forme
		$$E : \; y^2z = x^3+uxz^2+vz^3$$
		dans l'espace projectif $\mathbb{P}^2$ avec coordonn\'{e}es homog\`{e}nes $[x:y:z]$. Dans l'espace affine $\mathbb{A}^2\subset \mathbb{P}^2$ o\`{u} $z\neq 0$, cela correspond \`{a} l'\'{e}quation (obtenue en posant $z=1$ dans l'\'{e}quation de $E$) :
		$$E': \; y^2=x^3+ux+v.$$
		Le seul point de $E\setminus E'$ (\enquote{\`{a} l'infini}) est $[0:1:0]$. Si les param\`{e}tres $u,v\in\mathbb{Q}$ sont tels que $x^3+ux+v$ est \`{a} racines simples 
		on dit que $E$ est une \emph{courbe elliptique}. La figure suivante 
		montre $E(\RR)$ (le point \`{a} l'infini correspond \`{a} la direction verticale), et $E(\CC)$, une surface de Riemann compacte de genre $1$.
		\end{exb}
		
		\begin{figure}[h!]
		\includegraphics[scale=0.37]{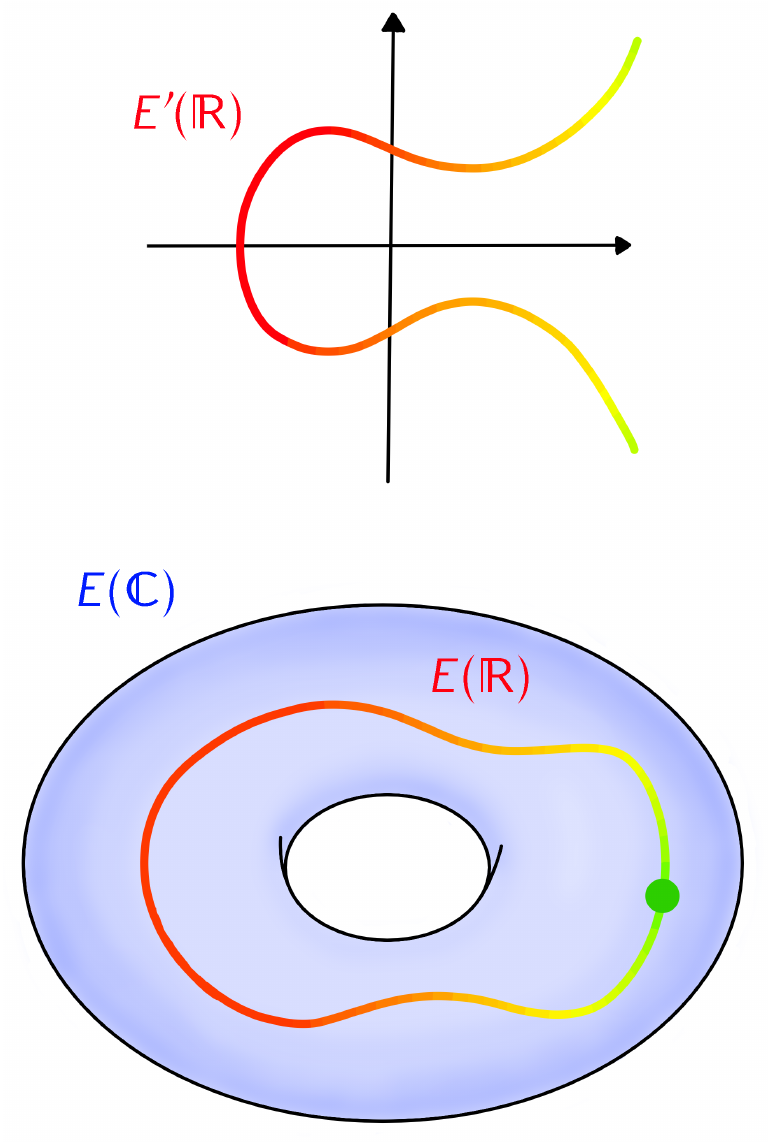}
		\end{figure}

	\subsection{Cohomologie des vari\'{e}t\'{e}s alg\'{e}briques}\label{par: cohomologie des varietes algebriques}
		
		Un outil a pris une importance consid\'{e}rable au \siecle{20} dans l'\'{e}tude (entre autres) des vari\'{e}t\'{e}s alg\'{e}briques : les \emph{groupes de cohomologie}, qui sont des espaces vectoriels de dimension finie
		$$\H^0(X), \H^1(X), \H^2(X), \ldots $$
		associ\'{e}s \`{a} une vari\'{e}t\'{e} alg\'{e}brique $X$, et qui encodent de mani\`{e}re concise certaines \enquote{obstructions} li\'{e}es \`{a} la g\'{e}om\'{e}trie de $X$. \`{A} la base de la th\'{e}orie des motifs figure le fait \'{e}tonnant qu'il y a plusieurs fa\c{c}ons de mettre en place un tel outil, dont les deux suivantes (voir plus bas pour des exemples de calculs).
		
		\begin{enumerate}
		\item[(a)] La \emph{cohomologie de Betti}, ou \emph{singuli\`{e}re}, disponible si $k$ est un sous-corps de $\mathbb{C}$, mise en place par Betti, Poincar\'{e}, E. Noether, et d'autres. Le groupe de cohomologie de Betti $\H^n_{\B}(X)$ est un $\QQ$-espace vectoriel de dimension finie, dual de l'homologie singuli\`{e}re \`{a} coefficients dans $\QQ$ de l'espace topologique $X(\CC)$ :
		$$\H^n_{\B}(X) = \H_n(X(\CC);\QQ)^\vee.$$
		Rappelons que l'homologie singuli\`{e}re est d\'{e}finie comme le quotient de l'espace des $n$-cycles topologiques par celui des $n$-bords.
		\item[(b)] La \emph{cohomologie de de Rham (alg\'{e}brique)}, d\'{e}finie par Grothendieck \cite{grothendieckderham} et qui est disponible si $k$ est de caract\'{e}ristique z\'{e}ro. Le groupe de cohomologie de de Rham $\H^n_{\dR}(X)$ est un $k$-espace vectoriel de dimension finie, analogue purement alg\'{e}brique de la cohomologie de de Rham en g\'{e}om\'{e}trie diff\'{e}rentielle (quotient de l'espace des $n$-formes diff\'{e}rentielles ferm\'{e}es par celui des formes exactes), les formes diff\'{e}rentielles $\mathcal{C}^\infty$ \'{e}tant remplac\'{e}es par les formes diff\'{e}rentielles alg\'{e}briques.
		\end{enumerate}	
	
		Si $k$ est un sous-corps de $\CC$, un th\'{e}or\`{e}me de Grothendieck, combin\'{e} au th\'{e}or\`{e}me de de Rham classique, donne un isomorphisme canonique entre ces deux th\'{e}ories cohomologiques, une fois les coefficients \'{e}tendus aux nombres complexes :
		\begin{equation}\label{eq: B dR comp}
		\H^n_{\dR}(X)\otimes_k\CC \xrightarrow{\,\sim\,} \H^n_{\B}(X)\otimes_\QQ\CC.
		\end{equation}
		Cet isomorphisme est issu de l'int\'{e}gration des formes diff\'{e}rentielles (repr\'{e}sentants de classes de cohomologie de de Rham) le long de cycles sur $X(\CC)$ (repr\'{e}sentants de classes d'homologie singuli\`{e}re).
		
		Notamment, les groupes de cohomologie de Betti et de de Rham ont la m\^{e}me dimension... ce qui est une \enquote{co\"{i}ncidence} \'{e}tonnante puisqu'ils mesurent des obstructions de natures tr\`{e}s diff\'{e}rentes !
		
		\begin{exa}
		Le groupe $\H^1_{\B/\dR}(\mathbb{A}^1\setminus \{0\})$ est de dimension $1$. Du c\^{o}t\'{e} Betti, il y a en effet une unique obstruction \`{a} pouvoir contracter un lacet $\gamma$ dans $\mathbb{C}^*$, mesur\'{e}e par l'\emph{indice} $\operatorname{Ind}(\gamma)$, le nombre de fois que $\gamma$ tourne autour de $0$. Du c\^{o}t\'{e} de Rham, il y a une unique obstruction \`{a} l'existence d'une primitive d'une fonction m\'{e}romorphe avec un seul p\^{o}le en $0$, mesur\'{e}e par le \emph{r\'{e}sidu} $\operatorname{Res}(f)$. 
		Ces obstructions, topologique et analytique, sont reli\'{e}es par la formule des r\'{e}sidus, qui est un pr\'{e}curseur de l'isomorphisme de comparaison \eqref{eq: B dR comp} :
		$$\int_\gamma f(z)\operatorname{d}\!z = 2\mathrm{i}\pi \, \operatorname{Ind}(\gamma) \, \operatorname{Res}(f).$$
		\end{exa}
		
		\begin{exb}
		La courbe elliptique $E$ et la courbe elliptique \'{e}point\'{e}e $E'$ ont les m\^{e}mes groupes de cohomologie $\H^1_{\B/\dR}(E)\simeq \H^1_{\B/\dR}(E')$, de dimension $2$. Une base de $\H^1_{\B}(E)^\vee = \H_1(E(\CC);\QQ)$ est donn\'{e}e par les classes de deux lacets $\alpha$ et $\beta$ comme dans la figure suivante. Une base de $\H^1_{\dR}(E')$ est donn\'{e}e\footnote{La cohomologie de de Rham alg\'{e}brique de $E'$ est plus facile \`{a} d\'{e}crire que celle de $E$ gr\^{a}ce au fait que $E'$ est affine~: toutes les classes de cohomologie peuvent \^{e}tre repr\'{e}sent\'{e}es par des formes alg\'{e}briques globales.} par les classes des formes diff\'{e}rentielles $\operatorname{d}\!x/y$ (qui s'\'{e}tend en une forme diff\'{e}rentielle sur $E$) 
		et $x\operatorname{d}\!x/y$ (qui a un p\^{o}le d'ordre $2$ au point \`{a} l'infini).
		\end{exb}
		
		\begin{figure}[h!]
		\includegraphics[scale=0.45]{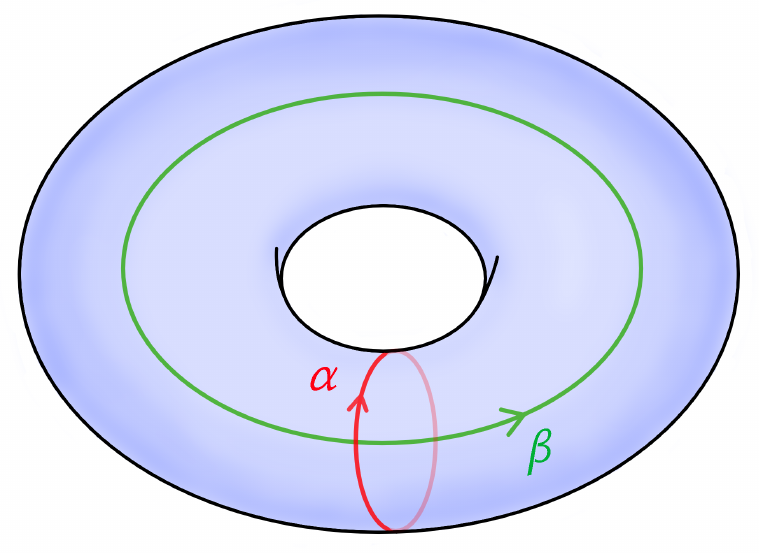}
		\end{figure}
		
		Plus g\'{e}n\'{e}ralement, pour une surface de Riemann compacte $X$ de genre $g$, le $\H^1_{\B/\dR}(X)$ est de dimension $2g$. Du c\^{o}t\'{e} Betti, un syst\`{e}me maximal de coupures de $X$ ne remplit que la \emph{moiti\'{e}} de ce groupe de cohomologie -- de m\^{e}me du c\^{o}t\'{e} de Rham avec une base des diff\'{e}rentielles de premi\`{e}re esp\`{e}ce. Ce ph\'{e}nom\`{e}ne est un cas particulier de la th\'{e}orie de Hodge, dont il sera question au \S\ref{par: hodge}.
		
		Les th\'{e}ories de Betti et de de Rham partagent des caract\'{e}ristiques formelles importantes, compatibles \`{a} l'isomorphisme de comparaison \eqref{eq: B dR comp} :
		\begin{itemize}
		\item la fonctorialit\'{e} (contravariante) : un morphisme de vari\'{e}t\'{e}s alg\'{e}briques $f:X\to Y$  induit une application lin\'{e}aire dans l'autre sens
		\begin{equation}\label{eq: fonctorialite cohomologie}
		f^*:\H^n(Y)\To \H^n(X)\ ;
		\end{equation}
		\item une structure d'alg\`{e}bre gradu\'{e}e commutative sur $\H^\bullet(X)$ induite par le \emph{cup-produit}
		$$\H^i(X)\otimes \H^j(X)\To \H^{i+j}(X) \ ;$$
		\item la \emph{formule de K\"{u}nneth} :
		$$\H^n(X\times Y) \simeq \bigoplus_{i+j=n} \H^i(X)\otimes \H^j(Y)\ ;$$
		\item la \emph{dualit\'{e} de Poincar\'{e}} : si $X$ est une vari\'{e}t\'{e} alg\'{e}brique projective lisse\footnote{Une vari\'{e}t\'{e} alg\'{e}brique projective est dite lisse si elle v\'{e}rifie le crit\`{e}re jacobien usuel de la g\'{e}om\'{e}trie diff\'{e}rentielle.} et connexe de dimension $d$ alors $\H^{2d}(X)$ est de dimension $1$ et l'accouplement 
		$$\H^i(X)\otimes \H^{2d-i}(X)\To \H^{2d}(X)$$
		donn\'{e} par le cup-produit est non d\'{e}g\'{e}n\'{e}r\'{e} pour tout entier naturel $i\in\{0,\ldots,2d\}$.
		\end{itemize}
		
		Mentionnons enfin, en anticipant le \S\ref{par: cycles algebriques} o\`{u} les d\'{e}finitions correspondantes seront donn\'{e}es, une caract\'{e}ristique essentielle de la cohomologie, qui la relie aux cycles alg\'{e}briques :
		\begin{itemize}
		\item des applications \emph{classe de cycle} issues des groupes de Chow, pour $X$ lisse :
		\begin{equation}\label{eq: classe de cycle general}
		\operatorname{CH}^r(X)_\QQ \To \H^{2r}(X) \quad , \quad Z\mapsto [Z].
		\end{equation}
		\end{itemize}
		
	\subsection{P\'{e}riodes}\label{par: periodes}
		
		 L'isomorphisme de comparaison \eqref{eq: B dR comp} peut \^{e}tre repr\'{e}sent\'{e} par une matrice, quitte \`{a} choisir une $k$-base de la cohomologie de de Rham et une $\QQ$-base de la cohomologie de Betti. Dans le cas o\`{u} $k$ est un corps de nombres (extension finie de $\QQ$) une telle matrice s'appelle \emph{matrice des p\'{e}riodes} et ses coefficients sont appel\'{e}s \emph{p\'{e}riodes}. Concr\`{e}tement, ce sont des int\'{e}grales de formes diff\'{e}rentielles alg\'{e}briques le long de cycles trac\'{e}s sur $X(\CC)$.
		
		\begin{exa}
		Dans le cas de $\H^1(\mathbb{A}^1\setminus \{0\})$, il y a essentiellement une seule p\'{e}riode, associ\'{e}e \`{a} la classe de la forme diff\'{e}rentielle $\frac{\operatorname{d}\!z}{z}$ et \`{a} la classe du lacet $\{|z|=1\}$ parcouru dans le sens trigonom\'{e}trique :
		$$\int_{|z|=1}\frac{\operatorname{d}\!z}{z} = 2\mathrm{i}\pi.$$
		\end{exa}
		
		\begin{exb}
		Une matrice des p\'{e}riodes de $\H^1(E)$ est de taille $2\times 2$.  Ses quatre coefficients sont des \enquote{int\'{e}grales elliptiques}, par exemple
		$$\int_a^{+\infty}\frac{\operatorname{d}\!x}{y} = \int_a^{+\infty} \frac{\operatorname{d}\!x}{\sqrt{x^3+ux+v}}$$
		o\`{u} $a$ est la plus grande racine r\'{e}elle de $x^3+ux+v$. Ce sont des p\'{e}riodes, au sens usuel, des fonctions elliptiques, variantes des fonctions trigonom\'{e}triques o\`{u} l'ellipse remplace le cercle. Des int\'{e}grales similaires donnent aussi la formule exacte pour la p\'{e}riode du pendule pesant en m\'{e}canique newtonienne.
		\end{exb}
		
		Les int\'{e}grales vues jusqu'\`{a} pr\'{e}sent se calculent sur des domaines sans bord. Il est naturel de consid\'{e}rer des int\'{e}grales plus g\'{e}n\'{e}rales en incluant les groupes de \emph{cohomologie relative} $\H^n(X,Y)$ o\`{u} $Y$ est une sous-vari\'{e}t\'{e} de $X$ (o\`{u} vit le bord du domaine d'int\'{e}gration). Par exemple, le groupe de cohomologie relative $\H^1(\mathbb{A}^1\setminus \{0\}, \{1,a\})$, pour un $a\in\QQ\setminus \{0,1\}$, contient dans sa matrice des p\'{e}riodes l'int\'{e}grale
		$$\int_1^a\frac{\operatorname{d}\!t}{t} = \log(a).$$
		On obtient alors une notion plus g\'{e}n\'{e}rale de p\'{e}riode, qui est \'{e}quivalente \`{a} la notion \'{e}l\'{e}mentaire suivante d\'{e}finie par Kontsevich et Zagier \cite{kontsevichzagier}.
		
		\begin{defi}
		Une \emph{p\'{e}riode} est un nombre complexe dont les parties r\'{e}elle et imaginaire s'\'{e}crivent sous la forme d'int\'{e}grales absolument convergentes
		$$\int_\sigma f(x_1,\ldots,x_n) \operatorname{d}\!x_1\cdots \operatorname{d}\!x_n$$
		o\`{u} $f\in\QQ(x_1,\ldots,x_n)$ et $\sigma\subset\mathbb{R}^n$ est un domaine semi-alg\'{e}brique d\'{e}fini sur $\QQ$, c'est-\`{a}-dire d\'{e}fini par un nombre fini d'in\'{e}galit\'{e}s $g\geq 0$ avec $g\in\QQ[x_1,\ldots,x_n]$.
		\end{defi}
		
		Les p\'{e}riodes forment un sous-anneau de $\mathbb{C}$ qui contient le corps $\overline{\QQ}$ des nombres alg\'{e}briques et qui est d\'{e}nombrable. On ne conna\^{i}t pourtant aucun nombre complexe \enquote{int\'{e}ressant} qui ne soit pas une p\'{e}riode, m\^{e}me s'il est conjectur\'{e} que la base $e$ du logarithme n\'{e}p\'{e}rien n'est pas une p\'{e}riode. En plus de $\pi$, des int\'{e}grales elliptiques, et des logarithmes, notons que la valeur de la fonction z\^{e}ta en un entier $n\geq 2$ est une p\'{e}riode, gr\^{a}ce \`{a} la formule
		$$\zeta(n) = \int_{[0,1]^n}\frac{\operatorname{d}\!x_1\cdots \operatorname{d}\!x_n}{1-x_1\cdots x_n}\cdot$$
		Elle appara\^{i}t dans la matrice des p\'{e}riodes d'un groupe de cohomologie relative associ\'{e} \`{a} un arrangement de sous-vari\'{e}t\'{e}s de $\mathbb{A}^n$.

		Par ailleurs, les p\'{e}riodes apparaissent naturellement en physique : les \emph{int\'{e}grales de Feynman}, qui calculent les probabilit\'{e}s des interactions entre particules \'{e}l\'{e}mentaires, s'expriment en fonction de p\'{e}riodes.
		
		L'\'{e}tude de l'arithm\'{e}tique des p\'{e}riodes est une des raisons d'\^{e}tre de la th\'{e}orie des motifs. 
		
	\subsection{Th\'{e}orie de Hodge}\label{par: hodge}
	
		Focalisons-nous un instant sur la cohomologie de Betti, c'est-\`{a}-dire sur les invariants topologiques des vari\'{e}t\'{e}s alg\'{e}briques complexes (d\'{e}finies sur $k=\CC$). Une question centrale est : comment la structure de vari\'{e}t\'{e} alg\'{e}brique contraint-elle la topologie ? La th\'{e}orie de Hodge est un ensemble d'outils pour r\'{e}pondre \`{a} cette question. Elle remonte au r\'{e}sultat fondateur suivant de Hodge \cite{hodge}, dont il est important de pr\'{e}ciser qu'il repose sur des m\'{e}thodes analytiques (g\'{e}om\'{e}trie k\"{a}hl\'{e}rienne).
		
		\begin{thm}[D\'{e}composition de Hodge]
		Soit $X$ une vari\'{e}t\'{e} alg\'{e}brique complexe projective et lisse, soit $n$ un entier naturel. On a une d\'{e}composition
		$$\H^n_{\B}(X)\otimes_{\QQ}\CC = \bigoplus_{\substack{p,q\geq 0\\ p+q=n}} \H^{p,q}(X)$$
		qui v\'{e}rifie la \enquote{sym\'{e}trie de Hodge}
		$$\overline{\H^{p,q}(X)} = \H^{q,p}(X).$$
		\end{thm}
		
		Concr\`{e}tement, l'espace $\H^{p,q}(X)$ est l'espace des classes de cohomologie qui peuvent \^{e}tre repr\'{e}sent\'{e}es par une forme diff\'{e}rentielle $\mathcal{C}^\infty$ complexe de type $(p,q)$, c'est-\`{a}-dire \'{e}crite en coordonn\'{e}es holomorphes locales $z_i$ avec $p$ occurrences de $\operatorname{d}\!z_i$ et $q$ occurrences de $\operatorname{d}\!\overline{z}_i$. La d\'{e}composition de Hodge est une structure suppl\'{e}mentaire contraignante sur la cohomologie de Betti ; par exemple, la sym\'{e}trie de Hodge implique que $\H^n_{\B}(X)$ est de dimension paire si $n$ est impair.
		
		Comme  $\H^{1,0}(X)$ est l'espace des $1$-formes holomorphes sur $X$, cela explique, dans le cas des surfaces de Riemann (vari\'{e}t\'{e}s alg\'{e}briques complexes projectives lisses de dimension $1$), que la dimension de $\H^1_{\B}(X)$ est le double du genre de $X$.
		
	\subsection{Cycles alg\'{e}briques}\label{par: cycles algebriques}
	
		La g\'{e}om\'{e}trie alg\'{e}brique regorge de probl\`{e}mes \'{e}num\'{e}ratifs 
		comme~: combien de droites de l'espace projectif de dimension $3$ sont contenues dans une surface cubique lisse donn\'{e}e ? (R\'{e}ponse : $27$) 
		
		Une approche moderne \`{a} ce genre de question est la th\'{e}orie de l'intersection, qui est l'\'{e}tude des cycles alg\'{e}briques. Un cycle alg\'{e}brique (ici, \`{a} coefficients dans $\QQ$) de codimension $r$ dans une vari\'{e}t\'{e} alg\'{e}brique $X$ 
		est une combinaison lin\'{e}aire formelle
		$$a_1Z_1+\cdots +a_nZ_n$$
		o\`{u} les $a_i\in\QQ$ et les $Z_i$ sont des sous-vari\'{e}t\'{e}s de codimension $r$ de $X$. On d\'{e}finit le $r$-i\`{e}me \emph{groupe de Chow} $\operatorname{CH}^r(X)_\QQ$ comme le quotient de l'espace des cycles alg\'{e}briques de codimension $r$ par la relation d'\'{e}quivalence rationnelle. Cette relation identifie deux cycles alg\'{e}briques qui sont reli\'{e}s par une famille de cycles alg\'{e}briques param\'{e}tr\'{e}s par $\mathbb{P}^1$. L'intersection des cycles alg\'{e}briques donne une structure d'anneau gradu\'{e} \`{a} $\operatorname{CH}^\bullet(X)_\QQ$,  not\'{e}e $(Z,Z')\mapsto Z\cdot Z'$. En effet, le quotient par l'\'{e}quivalence rationnelle permet de \enquote{d\'{e}placer} deux cycles alg\'{e}briques jusqu'\`{a} ce qu'ils s'intersectent dans la codimension attendue. Par exemple, la figure suivante repr\'{e}sente une surface quadrique lisse $Q$ de l'espace projectif $\mathbb{P}^3$ et une droite $D$ contenue dans $Q$. Comme $D$ est rationnellement \'{e}quivalente \`{a} une droite $D'$ qui intersecte $Q$ en deux points distincts $P_1$ et $P_2$, on a
		$$Q\cdot D = Q\cdot D' = P_1+P_2.$$
		
		\begin{figure}[h!]
		\includegraphics[scale=0.8]{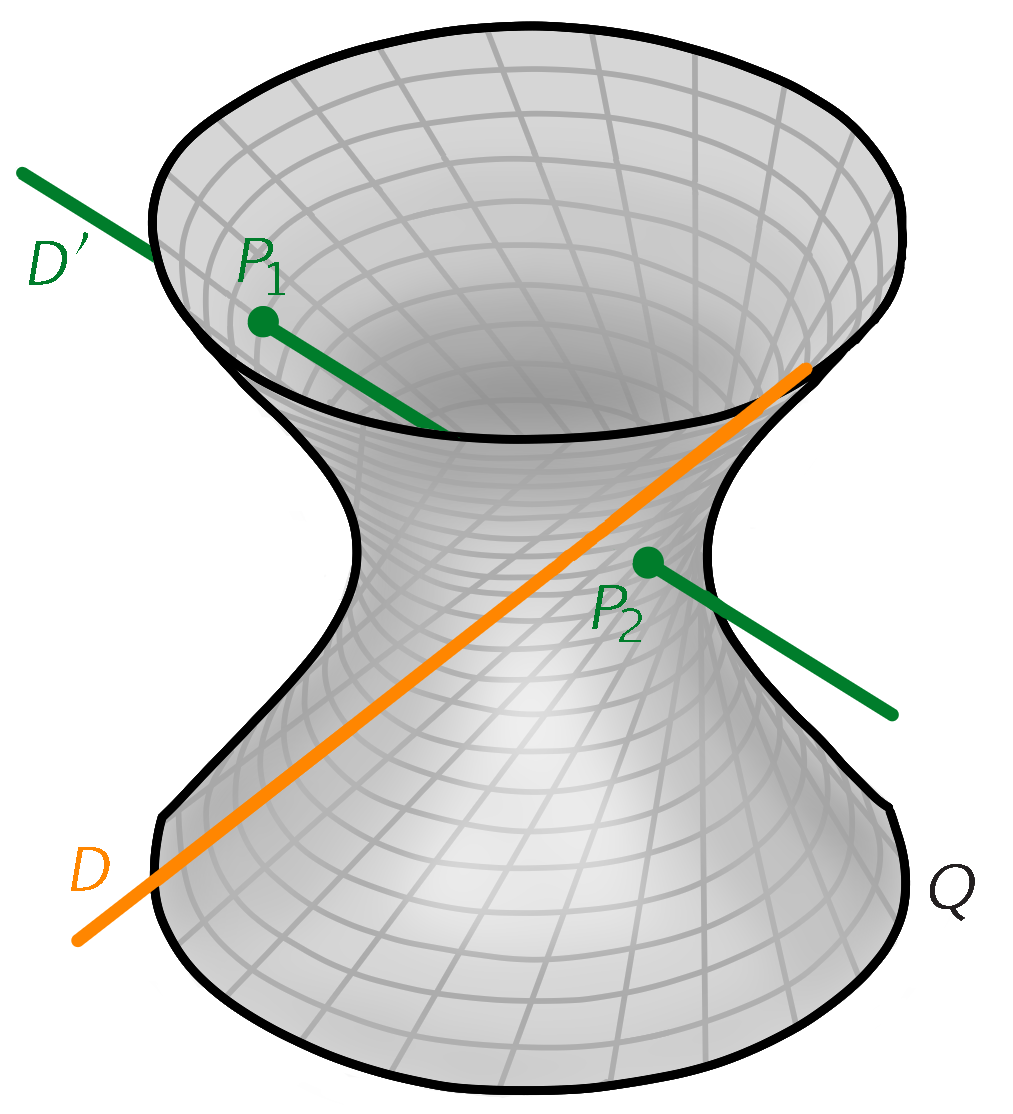}
		\end{figure}

		Le lien entre cycles alg\'{e}briques et cohomologie de Betti est donn\'{e} par l'application \emph{classe de cycle}, annonc\'{e}e au \S\ref{par: cohomologie des varietes algebriques}, pour $k$ un sous-corps de $\CC$ et $X$ une vari\'{e}t\'{e} lisse :
		\begin{equation}\label{eq: classe de cycle}
		\operatorname{CH}^r(X)_\QQ \To \H^{2r}_{\B}(X) \; , \; Z\mapsto [Z].
		\end{equation}
		Si $X$ est en plus projective  de dimension $d$ on a la dualit\'{e} de Poincar\'{e} $\H^{2r}_{\B}(X)\simeq \H_{2(d-r)}(X(\CC);\QQ)$, et $[Z]$ est repr\'{e}sent\'{e} par un cycle topologique de dimension $2(d-r)=\dim_{\mathbb{R}}(Z(\CC))$ d\'{e}fini par une triangulation quelconque de $Z(\CC)$.
		
		Dans le cas o\`{u} $X$ est projective et lisse, on voit facilement que l'image de l'application \eqref{eq: classe de cycle} est contenue dans la partie de type $(r,r)$ de la d\'{e}composition de Hodge. La conjecture de Hodge \cite{hodge} 
		est l'affirmation r\'{e}ciproque. \`{A} l'instar de l'hypoth\`{e}se de Riemann, il s'agit de l'un des sept \emph{Millenium Prize Problems} du Clay Mathematics Institute.
	
		\begin{conj}[Conjecture de Hodge]
		Soit $X$ une vari\'{e}t\'{e} alg\'{e}brique complexe projective et lisse, soit $r$ un entier naturel. Tout \'{e}l\'{e}ment de $\H^{2r}_{\B}(X)\cap \H^{r,r}(X)$ est la classe d'un cycle alg\'{e}brique de codimension $r$ dans $X$.
		\end{conj}

	\subsection{Les conjectures de Weil}\label{par: conjectures de weil}
		
		Les conjectures de Weil ont \'{e}t\'{e} une source incroyable de d\'{e}veloppements en g\'{e}om\'{e}trie arithm\'{e}tique. Elles ont trait \`{a} un probl\`{e}me tr\`{e}s ancien : compter le nombre de solutions enti\`{e}res d'\'{e}quations polynomiales \`{a} coefficients entiers (appel\'{e}es \emph{\'{e}quations diophantiennes}). Une strat\'{e}gie classique est de les r\'{e}duire modulo un nombre premier $p$ pour se ramener \`{a} des \'{e}quations dans le corps $\mathbb{F}_p=\ZZ/p\ZZ$ \`{a} $p$ \'{e}l\'{e}ments.  Rappelons que $\mathbb{F}_p$ a une unique extension de degr\'{e} $n$ pour chaque entier $n\geq 1$, qu'on note $\mathbb{F}_{p^n}$ et qui est l'unique corps \`{a} $p^n$ \'{e}l\'{e}ments.
		
		Adoptons un point de vue g\'{e}om\'{e}trique et consid\'{e}rons une vari\'{e}t\'{e} alg\'{e}brique $X$  d\'{e}finie sur $k=\mathbb{F}_q$, pour $q$ une puissance d'un nombre premier $p$.  Dans la lign\'{e}e des travaux de Riemann, et en suivant notamment E. Artin, Hasse, et Weil, on associe \`{a} $X$ sa \emph{fonction z\^{e}ta}, qui est la s\'{e}rie formelle
		$$Z(X,t) = \exp\left( \sum_{n\geq 1} \#(X(\mathbb{F}_{q^n}))\frac{t^n}{n}\right).$$
		On consid\'{e}rera \'{e}galement la variante
		$$\zeta_X(s) = Z(X,q^{-s}).$$
		dont les propri\'{e}t\'{e}s sont plus semblables \`{a} celles de la fonction z\^{e}ta de Riemann.
		
		Si $X$ est une courbe projective lisse, F. K. Schmidt \cite{schmidt} d\'{e}montre en 1931 la rationalit\'{e} de $Z(X,t)$ :
		$$Z(X,t) = \frac{P_X(t)}{(1-t)(1-qt)}\ ,$$
		o\`{u} $P_X(t)$ est un polyn\^{o}me. De plus, le degr\'{e} de $P_X(t)$ est $2g$ o\`{u} $g$ est le genre de $X$, d\'{e}fini comme la dimension de l'espace des $1$-formes alg\'{e}briques sans p\^{o}le sur $X$. Comme on le verra, l'apparition de ce nombre $2g$, qui est la dimension du $\H^1_{\B/\dR}$ d'une courbe projective lisse, n'est pas un hasard et a bien une explication cohomologique. Schmidt d\'{e}montre aussi que si l'on \'{e}crit
		$$P_X(t) = \prod_{i=1}^{2g}(1-\alpha_it)$$
		alors $\alpha\mapsto q/\alpha$ permute les $\alpha_i$. Ce r\'{e}sultat est un analogue de l'\'{e}quation fonctionnelle de la fonction z\^{e}ta de Riemann puisqu'il implique une relation entre $\zeta_X(s)$ et $\zeta_X(1-s)$.
		
		Dans le cas du genre $1$, Hasse \cite{hassecrelle} d\'{e}montre en 1933 que si l'on \'{e}crit $P_X(t)=(1-\alpha t)(1-\beta t)$, alors
		$$|\alpha|=|\beta|=\sqrt{q}.$$
		Ce r\'{e}sultat est un analogue de l'hypoth\`{e}se de Riemann: les racines de $\zeta_X$ sont sur la droite $\operatorname{Re}(s)=\frac{1}{2}$. Au cours des ann\'{e}es 1940, Weil \cite{weilnote} g\'{e}n\'{e}ralise le th\'{e}or\`{e}me de Hasse au cas du genre quelconque et en d\'{e}duit la \emph{borne de Hasse--Weil} sur le nombre de points d'une courbe projective lisse d\'{e}finie sur un corps fini :
		$$|\#(X(\mathbb{F}_q))-(q+1)| \leq 2g\sqrt{q}.$$
		
		Les conjectures de Weil \cite{weilBAMS}, \'{e}nonc\'{e}es en 1949, visent \`{a} g\'{e}n\'{e}raliser ces r\'{e}sultats \`{a} toutes les vari\'{e}t\'{e}s projectives et lisses d\'{e}finies sur un corps fini. Elles ont \'{e}t\'{e} d\'{e}montr\'{e}es l'une apr\`{e}s l'autre, par Dwork \cite{dwork} (rationalit\'{e} de $Z(X,t)$, 1959), Grothendieck et ses collaborateurs \cite{othersga4} (\'{e}quation fonctionnelle, lien avec la cohomologie de Betti, 1964), et enfin Deligne \cite{deligneweil1} (analogue de l'hypoth\`{e}se de Riemann, 1974). 
		
		Avant de donner une id\'{e}e de la preuve des conjectures de Weil, mentionnons la \emph{fonction z\^{e}ta globale} d'une vari\'{e}t\'{e} alg\'{e}brique $X$ d\'{e}finie sur $\QQ$, donn\'{e}e par le produit eul\'{e}rien
		$$\zeta_X(s) = \prod_{p \mbox{ \footnotesize{premier}}} Z(X_p,p^{-s}),$$
		o\`{u} $X_p$ est la r\'{e}duction de $X$ modulo $p$. La d\'{e}finition de $Z(X_p,t)$ doit \^{e}tre adapt\'{e}e pour un nombre fini de nombres premiers $p$, pour lesquels la r\'{e}duction de $X$ modulo $p$ n'est pas lisse. Si $X$ est un point, on retrouve la fonction z\^{e}ta de Riemann. Dans le cas d'une courbe elliptique $E$ on s'attend \`{a} ce que les propri\'{e}t\'{e}s analytiques de la fonction $\zeta_E$ refl\`{e}tent des propri\'{e}t\'{e}s arithm\'{e}tiques de $E$. Notamment, la conjecture de Birch et Swinnerton-Dyer -- un autre \emph{Millenium Prize Problem} -- relie le d\'{e}veloppement de Taylor de $\zeta_E(s)$ en $s=1$ \`{a} la structure du groupe des points rationnels $E(\QQ)$.
	
%

	\subsection{Encore plus de cohomologie}
	
		Weil \cite{weilICM} avait remarqu\'{e} que ses conjectures avaient un go\^{u}t cohomologique, par analogie avec le \emph{th\'{e}or\`{e}me du point fixe de Lefschetz}. Ce th\'{e}or\`{e}me compte, sous certaines hypoth\`{e}ses, le nombre de points fixes d'un endomorphisme $F$ d'un espace topologique $X$ comme somme altern\'{e}e des traces de l'endomorphisme $F_*$ induit par $F$ sur les diff\'{e}rents groupes d'homologie singuli\`{e}re de $X$ :

		$$\#(X^F) = \sum_{n\geq 0} (-1)^n \operatorname{Tr}(F_* :  \H_n(X;\QQ) \to \H_n(X;\QQ)).$$
		
		Quel rapport avec les conjectures de Weil ? Une vari\'{e}t\'{e} alg\'{e}brique $X$ d\'{e}finie sur $\mathbb{F}_q$ a un endomorphisme canonique $F$, l'\emph{endomorphisme de Frobenius}, qui \'{e}l\`{e}ve les coordonn\'{e}es \`{a} la puissance $q$, et $X(\mathbb{F}_{q^n})$ est l'ensemble des points fixes de $F^n$ dans l'ensemble $X(\overline{\mathbb{F}}_q)$, o\`{u} $\overline{\mathbb{F}}_q$ est une cl\^{o}ture alg\'{e}brique de $\mathbb{F}_q$. Weil remarque que si l'on disposait d'une th\'{e}orie cohomologique pour les vari\'{e}t\'{e}s sur $\mathbb{F}_q$ avec des propri\'{e}t\'{e}s formelles analogues \`{a} celles de la cohomologie singuli\`{e}re, on pourrait en d\'{e}duire une preuve de ses conjectures -- \`{a} l'exception notable de l'analogue de l'hypoth\`{e}se de Riemann. Notamment, le polyn\^{o}me $P_X(t)$ dont il a \'{e}t\'{e} question dans le \S\ref{par: conjectures de weil} s'interpr\'{e}terait, via un th\'{e}or\`{e}me du point fixe \`{a} la Lefschetz, comme le polyn\^{o}me caract\'{e}ristique de l'endomorphisme induit par le Frobenius sur $\H^1(X)$, et l'\'{e}quation fonctionnelle de la fonction $\zeta_X$ serait une cons\'{e}quence de la dualit\'{e} de Poincar\'{e}.
		
		Une telle th\'{e}orie de cohomologie a \'{e}t\'{e} d\'{e}velopp\'{e}e par Grothendieck et ses collaborateurs, notamment M. Artin, dans les ann\'{e}es 1960 \cite{othersga4}, et continue la liste des th\'{e}ories de cohomologie commenc\'{e}e au \S\ref{par: cohomologie des varietes algebriques} :
		\begin{enumerate}
		\item[(c)] La \emph{cohomologie \'{e}tale $\ell$-adique}, disponible pour tout corps $k$, qui produit des espaces vectoriels $\H_{\textnormal{\'{e}t},\ell}^n(X)$ sur le corps $\QQ_\ell$ des nombres $\ell$-adiques, o\`{u} $\ell$ est un nombre premier diff\'{e}rent de la caract\'{e}ristique de $k$.
		\end{enumerate}
		
		De mani\`{e}re analogue \`{a} la comparaison entre cohomologie de de Rham et de Betti \eqref{eq: B dR comp}, on a un isomorphisme de comparaison, d\^{u} \`{a} M. Artin, pour $k$ un sous-corps de $\CC$ :
		$$\H_{\textnormal{\'{e}t},\ell}^n(X) \stackrel{\sim}{\To} \H_{\B}^n(X)\otimes_\QQ \QQ_\ell.$$
		
		Par ailleurs, les groupes de cohomologie \'{e}tale $\ell$-adique ont une structure suppl\'{e}mentaire tr\`{e}s riche : ils sont naturellement des repr\'{e}sentations continues du groupe de Galois absolu du corps de d\'{e}finition $k$. Cela ouvre la voie \`{a} une \'{e}tude g\'{e}om\'{e}trique des groupes de Galois.
		
		\subsection{Vers les motifs}
		
			Betti, de Rham, \'{e}tale $\ell$-adique : trois th\'{e}ories de cohomologie\footnote{Mais aussi la \emph{cohomologie cristalline}, qui remplace la cohomologie de de Rham sur un corps parfait de caract\'{e}ristique positive.} pour les vari\'{e}t\'{e}s alg\'{e}briques qui \enquote{donnent les m\^{e}mes r\'{e}sultats\footnote{Un slogan \`{a} nuancer, les isomorphismes de comparaison cachant des ph\'{e}nom\`{e}ne subtils : Charles construit par exemple \cite{charlesconjugate} une vari\'{e}t\'{e} alg\'{e}brique $X$ d\'{e}finie sur un corps de nombres $k$ et deux plongements $k\subset \CC$ tels que les deux alg\`{e}bres de cohomologie de Betti $\H^\bullet_{\B}(X)$ correspondant aux deux plongements ne sont pas isomorphes.}}. Notamment, les groupes de cohomologie ont la m\^{e}me dimension et les m\^{e}mes propri\'{e}t\'{e}s formelles comme la dualit\'{e} de Poincar\'{e}.  Cependant, chaque th\'{e}orie a ses sp\'{e}cificit\'{e}s et ses structures suppl\'{e}mentaires, comme la d\'{e}composition de Hodge en cohomologie de Betti et l'action galoisienne en cohomologie \'{e}tale $\ell$-adique. De plus, les isomorphismes de comparaison entre th\'{e}ories de cohomologie font appara\^{i}tre une richesse arithm\'{e}tique, comme les p\'{e}riodes pour la comparaison entre Betti et de Rham \eqref{eq: B dR comp}.
			
			Cette profusion de th\'{e}ories de cohomologie est un luxe : on peut transf\'{e}rer intuition et techniques d'une th\'{e}orie \`{a} l'autre. Ce type de transfert est au c\oe ur de la strat\'{e}gie \'{e}nonc\'{e}e par Weil, et compl\'{e}t\'{e}e par Serre \cite{serreanalogues} et Grothendieck, pour attaquer ses conjectures. On verra aussi plus bas la notion de poids et la th\'{e}orie de Hodge mixte, d\'{e}velopp\'{e}e par Deligne \cite{delignehodge1} \`{a} partir de l'analogie entre cohomologie de Betti et cohomologie \'{e}tale $\ell$-adique.
		
			Il y a aussi des aspects plus frustrants de l'abondance des th\'{e}ories de cohomologie. Notamment, on a du mal \`{a} comparer les repr\'{e}sentations galoisiennes donn\'{e}es par la th\'{e}orie \'{e}tale $\ell$-adique pour diff\'{e}rents nombres premiers $\ell$. 
			
			Ces consid\'{e}rations ont amen\'{e} Grothendieck \`{a} d\'{e}gager la notion de motif comme th\'{e}orie de cohomologie \enquote{universelle}. Une intuition importante de Grothendieck est que la structure de cette th\'{e}orie devrait \^{e}tre contr\^{o}l\'{e}e par les cycles alg\'{e}briques. Les motifs forment donc un cadre pour aborder des questions sur l'interaction entre cycles alg\'{e}briques et cohomologie(s).
		
%
%
%

\section{La th\'{e}orie des motifs}

	Nous introduisons maintenant la th\'{e}orie des motifs \`{a} proprement parler. Le lecteur impatient pourra  parcourir le \S\ref{par: morphismes entre motifs} avant de passer aux aspects galoisiens et aux applications aux p\'{e}riodes du \S\ref{par: aspects galoisiens}.

	\subsection{Motifs... et morphismes entre motifs}\label{par: morphismes entre motifs}
	
		Le \emph{motif} d'une vari\'{e}t\'{e} alg\'{e}brique $X$ est pens\'{e} par Grothendieck comme une version \enquote{universelle} de sa cohomologie, c'est-\`{a}-dire un objet $\M(X)$ qui contr\^{o}le les propri\'{e}t\'{e}s des groupes $\H^n(X)$ dans les diff\'{e}rentes th\'{e}ories de cohomologie. 
		
		Pourquoi d\'{e}sirer une telle notion alors qu'on pourrait se contenter de la cohomologie ? Une des raisons est qu'il arrive souvent que la cohomologie de deux vari\'{e}t\'{e}s alg\'{e}briques diff\'{e}rentes $X$ et $Y$ contienne un \enquote{morceau commun} qui est \enquote{le m\^{e}me} dans chaque th\'{e}orie : m\^{e}me dimension, mais aussi m\^{e}mes p\'{e}riodes, m\^{e}me d\'{e}composition de Hodge, m\^{e}me action du groupe de Galois absolu, etc. Cela se traduit notamment par un facteur commun dans les fonctions $L$ associ\'{e}es \`{a} $X$ et $Y$. On a alors envie d'exprimer, et si possible d'expliquer, cette \enquote{co\"{i}ncidence} \`{a} un niveau conceptuel : le morceau commun est un objet en lui-m\^{e}me, un motif au sens math\'{e}matique, qui est un sous-objet commun aux objets $\M(X)$ et $\M(Y)$.
		
		L'exemple le plus simple de ce ph\'{e}nom\`{e}ne est que le $\H^0$ d'une vari\'{e}t\'{e} alg\'{e}brique connexe est toujours \enquote{le m\^{e}me}. L'explication est ici \'{e}l\'{e}mentaire : le morphisme de $X$ vers un point induit par fonctorialit\'{e} une application lin\'{e}aire $\H^0(\operatorname{point})\to \H^0(X)$ qui est un isomorphisme dans chaque th\'{e}orie de cohomologie. En g\'{e}n\'{e}ral, trouver des explications \`{a} de telles \enquote{co\"{i}ncidences} est beaucoup plus compliqu\'{e}. C'est un des buts de la th\'{e}orie des motifs.
		

		
		 On est donc amen\'{e} \`{a} donner un sens \`{a} la notion de \enquote{morceau commun} via une notion d'\emph{isomorphisme entre motifs}, et plus g\'{e}n\'{e}ralement \`{a} relier des motifs diff\'{e}rents via une notion de \emph{morphisme entre motifs}. La question 
		\begin{center} \enquote{Qu'est-ce qu'un motif ?} \end{center} 
		n'est donc pas pertinente sans la question
		\begin{center} \enquote{Qu'est-ce qu'un morphisme entre motifs ?}. \end{center}
		Dit autrement, on cherche une \emph{cat\'{e}gorie de motifs}. Rappelons qu'une \emph{cat\'{e}gorie} $\mathsf{C}$ est la donn\'{e}e 
		\begin{itemize}
		\item d'une classe d'objets (par exemple les ensembles, les espaces vectoriels, les groupes)
		\item et d'une notion de morphisme (par exemple les applications, les applications lin\'{e}aires, les morphismes de groupes), c'est-\`{a}-dire, d'ensembles de morphismes $\operatorname{Hom}_{\mathsf{C}}(M,N)$ entre deux objets $M$ et $N$, et d'une notion de composition des morphismes.
		\end{itemize}
		Il existe plusieurs approches des motifs, qui ne sont pas encore toutes unifi\'{e}es, et sont plus ou moins utiles en fonction du contexte. Listons quelques-unes des caract\'{e}ristiques attendues de \enquote{la} th\'{e}orie des motifs. 
		\begin{itemize}
		\item[(M1)] Les motifs sont les objets d'une cat\'{e}gorie $\mathsf{Mot}(k)$, la \emph{cat\'{e}gorie des motifs} sur $k$, qui est $\QQ$-lin\'{e}aire : les ensembles de morphismes sont des $\QQ$-espaces vectoriels et la composition des morphismes est bilin\'{e}aire.\footnote{Pour des motifs plus g\'{e}n\'{e}raux, l'anneau des coefficients n'est pas n\'{e}cessairement $\QQ$. On peut aussi vouloir remplacer le corps de d\'{e}finition $k$ par un anneau commutatif g\'{e}n\'{e}ral, voire une vari\'{e}t\'{e} alg\'{e}brique de param\`{e}tres.} Cet axiome est raisonnable puisque les motifs sont des abstractions des groupes de cohomologie, et donc des objets \enquote{lin\'{e}aires}.
		\item[(M2)] Pour chaque th\'{e}orie de cohomologie (Betti, de Rham, \'{e}tale $\ell$-adique) \`{a} coefficients dans un corps $F$ on a un \emph{foncteur de r\'{e}alisation}\footnote{Un foncteur $F$ d'une cat\'{e}gorie $\mathsf{C}$ vers une cat\'{e}gorie $\mathsf{D}$ associe \`{a} tout objet $M$ de $\mathsf{C}$ un objet $F(M)$ de $\mathsf{D}$, et \`{a} tout morphisme $M\to N$ dans $\mathsf{C}$ un morphisme $F(M)\to F(N)$ dans $\mathsf{D}$, de mani\`{e}re compatible \`{a} la composition des morphismes.}
		$$\mathsf{Mot}(k) \longrightarrow \mathsf{Vect}_F$$
		vers la cat\'{e}gorie des espaces vectoriels sur $F$ de dimension finie. 
		\item[(M3)] Pour toute vari\'{e}t\'{e} alg\'{e}brique $X$ d\'{e}finie sur $k$ et tout entier naturel $n$, on a un objet $\operatorname{M}^n(X)$ de $\mathsf{Mot}(k)$, le \emph{motif} de $X$ en degr\'{e} $n$, dont les r\'{e}alisations (images par un des foncteurs de r\'{e}alisation) sont les groupes de cohomologie $\H^n(X)$ dans les diff\'{e}rentes th\'{e}ories. (Plus g\'{e}n\'{e}ralement on veut aussi avoir \`{a} disposition des objets $\M^n(X,Y)$ jouant le r\^{o}le de groupes de cohomologie relative, m\^{e}me si on n'insistera pas sur cet aspect.) Plus formellement, on a des foncteurs contravariants\footnote{Un foncteur contravariant renverse le sens des fl\`{e}ches : \`{a} un morphisme $M\to N$ il associe un morphisme $F(N)\to F(M)$ \enquote{dans l'autre sens}, de mani\`{e}re compatible \`{a} la composition des morphismes.}
		\begin{equation}\label{eq: foncteur Var Mot}
		\mathsf{Var}(k) \longrightarrow \mathsf{Mot}(k) \;\; , \;\; X\mapsto \operatorname{M}^n(X)
		\end{equation}
		qui factorisent les foncteurs de cohomologie. 
		On peut \^{e}tre plus ambitieux et demander que la cat\'{e}gorie $\mathsf{Mot}(k)$ soit la cat\'{e}gorie \emph{universelle} qui factorise tous les foncteurs de cohomologie d'une classe donn\'{e}e.
		\item[(M4)] La cat\'{e}gorie $\mathsf{Mot}(k)$ a une structure de cat\'{e}gorie mono\"{i}dale sym\'{e}trique (c'est-\`{a}-dire une notion de produit tensoriel) et les motifs $\operatorname{M}^n(X)$ satisfont \`{a} une formule de K\"{u}nneth, et \`{a} la dualit\'{e} de Poincar\'{e} si $X$ est projective et lisse.
		\end{itemize}
		
		Enfin, ajoutons une condition qui n'est pas impos\'{e}e dans certaines approches de la th\'{e}orie des motifs, et que nous d\'{e}velopperons au \S\ref{par: aspects galoisiens} avec des applications aux p\'{e}riodes.
		\begin{itemize}
		\item[(M5)] La cat\'{e}gorie $\mathsf{Mot}(k)$ est tannakienne, avec les foncteurs de r\'{e}alisation pour foncteurs fibre.
		\end{itemize}
		Contentons-nous de noter qu'une cat\'{e}gorie tannakienne est notamment ab\'{e}lienne (existence de noyaux et images avec les propri\'{e}t\'{e}s attendues).

	\subsection{Motifs et cycles alg\'{e}briques}
	
		
		La caract\'{e}ristique cruciale des motifs est l'axiome suivant.
		
		\begin{itemize}
		\item[(M6)] Les morphismes dans la cat\'{e}gorie $\mathsf{Mot}(k)$ sont reli\'{e}s aux cycles alg\'{e}briques.
		\end{itemize}
		
		Pr\'{e}cisons ce desideratum. Les applications lin\'{e}aires \eqref{eq: fonctorialite cohomologie} ne sont pas les seules disponibles dans toutes les th\'{e}ories cohomologiques. En effet, les applications classe de cycle \eqref{eq: classe de cycle general} permettent de d\'{e}finir des applications lin\'{e}aires entre groupes de cohomologie gr\^{a}ce \`{a} la notion de correspondance.  Une \emph{correspondance} entre deux vari\'{e}t\'{e}s alg\'{e}briques $X$ et $Y$ est un cycle alg\'{e}brique dans le produit $X\times Y$. Notons $\operatorname{pr}_X:X\times Y\to X$ et $\operatorname{pr}_Y:X\times Y\to Y$ les deux projections. Si $Z$ est une correspondance de codimension $r$ entre $X$ et $Y$ projectives lisses, on peut consid\'{e}rer la composition :
		\begin{equation}\label{eq: morphisme corr}
		\H^n(Y) \xrightarrow{(\operatorname{pr}_Y)^*}   \H^n(X\times Y)  \xrightarrow{[Z]\cdot(-)}  \H^{n+2r}(X\times Y)   \xrightarrow{(\operatorname{pr}_X)_*}   \H^{n+2r-2\dim(Y)}(X)  
		\end{equation}
		o\`{u} $(\operatorname{pr}_X)_*$ est l'\enquote{int\'{e}gration le long des fibres}, issue de la dualit\'{e} de Poincar\'{e}. Le graphe d'un morphisme $f:X\to Y$ est un cas particulier de correspondance de codimension $r=\dim(Y)$ entre $X$ et $Y$, pour lequel \eqref{eq: morphisme corr} est l'application lin\'{e}aire \eqref{eq: fonctorialite cohomologie} induite par $f$.
		
		Une intuition importante de Grothendieck est que, dans le cadre des vari\'{e}t\'{e}s projectives lisses, les applications lin\'{e}aires \eqref{eq: morphisme corr} sont \emph{les seules} qui doivent \^{e}tre consid\'{e}r\'{e}s comme communes \`{a} toutes les th\'{e}ories cohomologiques, et sont donc la \emph{d\'{e}finition} des morphismes entre motifs de vari\'{e}t\'{e}s projectives lisses. (Nous verrons que le cas des vari\'{e}t\'{e}s g\'{e}n\'{e}rales est plus subtil.)

	\subsection{Motifs purs \`{a} la Grothendieck}
	
		La premi\`{e}re d\'{e}finition d'une cat\'{e}gorie des motifs est due \`{a} Grothendieck. Il s'agit d'une th\'{e}orie partielle puisqu'elle ne concerne que les vari\'{e}t\'{e}s alg\'{e}briques projectives et lisses. Les motifs correspondants sont appel\'{e}s \emph{purs} pour les distinguer des motifs g\'{e}n\'{e}raux qu'on appelle parfois \emph{mixtes}. On doit \`{a} Manin \cite{manincorrespondences} le premier texte sur les id\'{e}es de Grothendieck au sujet des motifs purs.
		
		Nous pr\'{e}sentons maintenant la construction (en trois \'{e}tapes) d'une des versions de la cat\'{e}gorie des motifs purs de Grothendieck. Il s'agit de la cat\'{e}gorie des \emph{motifs num\'{e}riques}, bas\'{e}e sur l'\'{e}quivalence num\'{e}rique entre cycles alg\'{e}briques. Pour une vari\'{e}t\'{e} alg\'{e}brique projective lisse $X$ de dimension $d$ et deux \'{e}l\'{e}ments $Z,Z'\in \operatorname{CH}^r(X)_\QQ$, on dit que $Z$ et $Z'$ sont \emph{num\'{e}riquement \'{e}quivalents}, et on \'{e}crit $Z\sim_{\operatorname{num}}Z'$, si pour tout $W\in \operatorname{CH}^{d-r}(X)_\QQ$ de codimension compl\'{e}mentaire, les nombres d'intersection $\#(Z\cdot W)$ et $\#(Z'\cdot W)$ co\"{i}ncident.
		
		\subsubsection{\'{E}tape 1}
		
		On consid\`{e}re la cat\'{e}gorie $\mathsf{C}$ dont les objets sont les vari\'{e}t\'{e}s alg\'{e}briques projectives lisses, o\`{u} l'on note $\operatorname{M}(X)$ l'objet correspondant \`{a} une vari\'{e}t\'{e} $X$ ; les morphismes sont donn\'{e}s par les correspondances modulo \'{e}quivalence num\'{e}rique :
		$$\operatorname{Hom}_{\mathsf{C}}(\operatorname{M}(Y),\operatorname{M}(X)) = \operatorname{CH}^{\dim(Y)}(X\times Y)_\QQ \,/ \sim_{\operatorname{num}} .$$
		La composition de correspondances $Z_{12}$ entre $X_1$ et $X_2$ et $Z_{23}$ entre $X_2$ et $X_3$ se fait par la formule
		$$Z_{12}\circ Z_{23} = (\operatorname{pr}_{13})_*((\operatorname{pr}_{12})^*(Z_{12})\cdot (\operatorname{pr}_{23})^*(Z_{23}) ),$$
		o\`{u} les $\operatorname{pr}_{ij}:X_1\times X_2\times X_3\to X_i\times X_j$ sont les projections. La cat\'{e}gorie $\mathsf{C}$ est $\QQ$-lin\'{e}aire et mono\"{i}dale sym\'{e}trique, le produit tensoriel \'{e}tant donn\'{e} par le produit des vari\'{e}t\'{e}s. L'objet correspondant \`{a} un point est not\'{e} $\QQ$, c'est l'unit\'{e} du produit tensoriel.
		
		On a un foncteur contravariant $\mathsf{Var}(k)\to \mathsf{C}$ qui envoie $X$ sur $\operatorname{M}(X)$ et un morphisme $f:X\to Y$ vers la classe du graphe de $f$.
	
		\subsubsection{\'{E}tape 2}
		
		L'objet $\M(X)$ de $\mathsf{C}$ joue le r\^{o}le de la somme directe des $\H^n(X)$ pour tous les degr\'{e}s $n$. On souhaiterait \enquote{d\'{e}couper} cet objet en une somme directe d'objets $\M^n(X)$. Plus g\'{e}n\'{e}ralement, on souhaiterait aussi \enquote{d\'{e}couper} chaque $\operatorname{M}^n(X)$ en somme directe d'objets plus petits (\enquote{particules \'{e}l\'{e}mentaires} ?) quand une telle d\'{e}composition existe dans chaque th\'{e}orie de cohomologie.
		
		On d\'{e}finit donc la cat\'{e}gorie $\mathsf{D}$ comme \'{e}tant la \emph{compl\'{e}tion pseudo-ab\'{e}lienne} de $\mathsf{C}$. Cela revient \`{a} rajouter formellement des objets qui jouent le r\^{o}le de noyaux ou images des projecteurs : pour tout objet $M$ de $\mathsf{C}$ et tout morphisme $e:M\to M$ qui v\'{e}rifie $e\circ e=e$, on a des objets $\ker(e)$ et $\mathsf{Im}(e)$ dans $\mathsf{D}$ et une d\'{e}composition $M=\ker(e)\oplus \operatorname{Im}(e)$.
		
		On obtient alors une d\'{e}composition $\operatorname{M}(\mathbb{P}^1)=\operatorname{M}^0(\mathbb{P}^1)\oplus \operatorname{M}^2(\mathbb{P}^1)$ dans $\mathsf{D}$, gr\^{a}ce au projecteur donn\'{e} par le cycle alg\'{e}brique $\mathbb{P}^1\times\{0\}\subset \mathbb{P}^1\times \mathbb{P}^1$. On montre facilement que $\operatorname{M}^0(\mathbb{P}^1)=\QQ$ est le motif trivial. M\^{e}me si $\H^2(\mathbb{P}^1)$ est de dimension $1$, le motif $\operatorname{M}^2(\mathbb{P}^1)$  n'est pas isomorphe au motif trivial. Une mani\`{e}re de le voir est que parmi les p\'{e}riodes de $\H^2(\mathbb{P}^1)$ figure $2\mathrm{i}\pi$, alors que les p\'{e}riodes du $\H^0$ du point sont des nombres rationnels. Ou encore, en th\'{e}orie de Hodge, $\H^2(\mathbb{P}^1)$ est de type $(1,1)$ alors que le $\H^0$ du point est de type $(0,0)$. On note
		$$\operatorname{M}^2(\mathbb{P}^1)=\QQ(-1),$$
		qu'on appelle le \emph{motif de Lefschetz}.
		
		\subsubsection{\'{E}tape 3}
		
		La cat\'{e}gorie des motifs num\'{e}riques, not\'{e}e $\mathsf{NumMot}(k)$, est d\'{e}finie \`{a} partir de $\mathsf{D}$ en rajoutant formellement un inverse tensoriel \`{a} $\QQ(-1)$, c'est-\`{a}-dire un objet $\QQ(1)$, appel\'{e} \emph{motif de Tate}, qui v\'{e}rifie
		$$\QQ(-1)\otimes \QQ(1) \simeq \QQ.$$
		Cette \'{e}tape a pour but d'induire sur $\mathsf{NumMot}(k)$ une notion de dualit\'{e} compatible \`{a} la dualit\'{e} de Poincar\'{e} en cohomologie : pour une vari\'{e}t\'{e} projective lisse $X$ de dimension $d$, on a 
		$$\M(X)^\vee = \M(X)\otimes \QQ(1)^{\otimes d} $$
		

			Le th\'{e}or\`{e}me suivant, conjectur\'{e} par Grothendieck et d\'{e}montr\'{e} par Jannsen \cite{jannsen} en 1992, justifie la construction.
		
			\begin{thm}
			La cat\'{e}gorie $\mathsf{NumMot}(k)$ est ab\'{e}lienne semi-simple\footnote{Un objet $M$ d'une cat\'{e}gorie ab\'{e}lienne $\mathsf{C}$ est dit simple si ses seuls sous-objets sont $0$ et $M$. On dit que $\mathsf{C}$ est semi-simple si tout objet de $\mathsf{C}$ est somme directe d'objets simples.}.
			\end{thm}
			
			La cat\'{e}gorie des motifs num\'{e}riques est donc une bonne candidate pour \enquote{la} cat\'{e}gorie des motifs purs. Comme on le verra plus bas, on ne sait cependant pas d\'{e}montrer qu'elle satisfait tous les axiomes (M1)-(M6).

		\subsubsection{Motifs de Chow, motifs homologiques}
		
			On peut r\'{e}p\'{e}ter la construction ci-dessus avec les variantes suivantes :
			
			\begin{itemize} 
			\item travailler avec les groupes de Chow non quotient\'{e}s, ce qui produit la cat\'{e}gorie des \emph{motifs de Chow} $\mathsf{CHMot}(k)$ ;
			\item remplacer l'\'{e}quivalence num\'{e}rique par l'\'{e}quivalence homologique $\sim_{\operatorname{hom}}$ (deux cycles alg\'{e}briques sont homologiquement \'{e}quivalents si leurs classes dans une th\'{e}orie cohomologique donn\'{e}e co\"{i}ncident), ce qui produit la cat\'{e}gorie des \emph{motifs homologiques} $\mathsf{HomMot}(k)$.
			\end{itemize}
			
			Des cycles alg\'{e}briques homologiquement \'{e}quivalents sont num\'{e}riquement \'{e}quivalents, d'o\`{u} des foncteurs
			$$\mathsf{CHMot}(k) \To \mathsf{HomMot}(k) \To \mathsf{NumMot(k)}.$$

		\subsection{Les conjectures standard}
		
			On souhaiterait d\'{e}finir, pour une th\'{e}orie de cohomologie donn\'{e}e \`{a} coefficients dans un corps $F$, un foncteur de r\'{e}alisation $\mathsf{NumMot}(k)\to \mathsf{Vect}_F$ par la formule 
			$$\operatorname{M}(X)\;\; \mapsto \;\; \H^\bullet(X)=\bigoplus_{n\geq 0}\H^n(X),$$
			et par \eqref{eq: morphisme corr} pour les morphismes. Mais apr\`{e}s avoir quotient\'{e} les groupes de Chow par l'\'{e}quivalence num\'{e}rique, cela n'a de sens que si deux cycles alg\'{e}briques num\'{e}riquement \'{e}quivalents sont homologiquement \'{e}quivalents, c'est-\`{a}-dire si la conjecture suivante est v\'{e}rifi\'{e}e.
		
			\begin{conj}
			Les relations d'\'{e}quivalence homologique et num\'{e}rique co\"{i}ncident : $\sim_{\operatorname{hom}} = \sim_{\operatorname{num}}$.
			\end{conj}
			
			On s'attend donc \`{a} une \'{e}quivalence de cat\'{e}gories : $\mathsf{HomMot}(k)\stackrel{?}{\simeq} \mathsf{NumMot}(k)$.  Notons que la cat\'{e}gorie $\mathsf{CHMot}(k)$ n'est pas ab\'{e}lienne en g\'{e}n\'{e}ral (mais on conjecture qu'elle l'est dans le cas d'un corps fini $k$). Son r\^{o}le dans la th\'{e}orie est clarifi\'{e} par la notion de motif mixte, qu'on d\'{e}veloppera au \S\ref{par: motifs mixtes}.
			
			Si l'on travaille avec les motifs homologiques pour assurer l'existence de foncteurs de r\'{e}alisation, un autre probl\`{e}me se pose (en plus du fait qu'on ne sait pas prouver que la cat\'{e}gorie $\mathsf{HomMot}(k)$ est ab\'{e}lienne). En effet, on ne sait pas en g\'{e}n\'{e}ral \enquote{d\'{e}couper} l'objet $\operatorname{M}(X)$ de $\mathsf{HomMot}(k)$ en une somme directe d'objets $\operatorname{M}^n(X)$ qui se r\'{e}aliseraient en les groupes de cohomologie $\H^n(X)$. Cela revient \`{a} la conjecture suivante, dont on peut montrer qu'elle est en fait une cons\'{e}quence de la pr\'{e}c\'{e}dente.

		\begin{conj}
		Soit $X$ une vari\'{e}t\'{e} alg\'{e}brique projective et lisse de dimension $d$, soit un entier $i\in\{0,\ldots,2d\}$. Le \enquote{projecteur de K\"{u}nneth}
		$$\H^\bullet(X) \twoheadrightarrow \H^i(X) \hookrightarrow \H^\bullet(X)$$
		est induit par une correspondance, c'est-\`{a}-dire un \'{e}l\'{e}ment de $\operatorname{CH}^d(X\times X)_\QQ$, via la construction \eqref{eq: morphisme corr}.
		\end{conj}
		
		Dans le cas o\`{u} $k$ est un corps fini, cette conjecture a \'{e}t\'{e} d\'{e}montr\'{e}e par Katz et Messing \cite{katzmessing}.
		
		Les deux conjectures pr\'{e}c\'{e}dentes font partie d'un ensemble coh\'{e}rent de \enquote{conjectures standard sur les cycles alg\'{e}briques} \'{e}nonc\'{e}es par Grothendieck \cite{grothendieckstandard, kleimanstandard} en 1969. Grothendieck voyait un double int\'{e}r\^{e}t \`{a} ces conjectures : elles permettaient de fonder la th\'{e}orie des motifs, et avaient comme cons\'{e}quence formelle la plus difficile des conjectures de Weil, l'analogue de l'hypoth\`{e}se de Riemann. Celle-ci d\'{e}coulerait en effet de r\'{e}sultats de positivit\'{e} de formes quadratiques sur les cycles alg\'{e}briques (modulo \'{e}quivalence), analogues du th\'{e}or\`{e}me de l'indice de Hodge en g\'{e}om\'{e}trie complexe.
		
		Malheureusement, tr\`{e}s peu de progr\`{e}s ont \'{e}t\'{e} faits sur les conjectures standard depuis lors~: prouver ces conjectures, ou encore la conjecture de Hodge, n\'{e}cessite de produire des cycles alg\'{e}briques, t\^{a}che pour laquelle on ne dispose pas de techniques g\'{e}n\'{e}rales. De plus, Deligne \cite{deligneweil1} a conclu la preuve des conjectures de Weil sans recourir \`{a} la strat\'{e}gie sugg\'{e}r\'{e}e par Grothendieck. Ironiquement, les travaux de Deligne sur les conjectures de Weil ont permis de faire des progr\`{e}s sur les conjectures standard, dont le th\'{e}or\`{e}me de Katz--Messing \'{e}voqu\'{e} ci-dessus. 
		
		Dans le dernier demi-si\`{e}cle, la th\'{e}orie des motifs s'est d\'{e}velopp\'{e}e dans une direction qui n'\'{e}tait pas celle initialement envisag\'{e}e par Grothendieck, contournant les conjectures standard. Celles-ci restent cependant un guide pr\'{e}cieux et une source importante de questions sur les cycles alg\'{e}briques, et les cat\'{e}gories de motifs purs \`{a} la Grothendieck sont un outil central pour y r\'{e}pondre. 
		
		Mentionnons au passage la th\'{e}orie des \emph{cycles motiv\'{e}s} d'Andr\'{e} \cite{andreinconditionnelle}, faisant suite \`{a} la notion de \emph{cycle de Hodge absolu} de Deligne \cite{delignehodgecycles}, et qui permet de mettre au point une th\'{e}orie inconditionnelle des motifs purs ayant les caract\'{e}ristiques attendues (M1)-(M6). L'id\'{e}e est de remplacer les cycles alg\'{e}briques par certaines classes de cohomologie qui viennent conjecturalement des cycles alg\'{e}briques.

	\subsection{Extensions entre motifs}
		
			Un pr\'{e}curseur de la notion de motif mixte est la notion de \emph{poids}, d\'{e}velopp\'{e}e par Deligne \cite{delignehodge1} d'apr\`{e}s des id\'{e}es de Grothendieck. L'id\'{e}e est que la cohomologie des vari\'{e}t\'{e}s alg\'{e}briques (dans toutes les th\'{e}ories) est munie d'une filtration canonique, la \emph{filtration par le poids}
			$$\cdots \subset W_{i-1} \H^n(X) \subset W_i\H^n(X)\subset \cdots \subset \H^n(X)$$		
			telle que le quotient $W_i/W_{i-1}$ est \enquote{pur de poids $i$}, c'est-\`{a}-dire s'exprime en fonction de $\H^i$ de vari\'{e}t\'{e}s alg\'{e}briques projectives lisses. Les groupes de cohomologie des vari\'{e}t\'{e}s alg\'{e}briques projectives lisses sont donc les \enquote{briques \'{e}l\'{e}mentaires} permettant de construire les groupes de cohomologie de toutes les vari\'{e}t\'{e}s alg\'{e}briques. En cohomologie \'{e}tale $\ell$-adique la filtration par le poids est induite par les valeurs propres de l'endomorphisme de Frobenius, alors qu'en cohomologie de Betti elle est au c\oe ur de la th\'{e}orie de Hodge mixte.
			
			\`{A} titre d'exemple, consid\'{e}rons une courbe alg\'{e}brique projective lisse $X$ et deux points rationnels $a\neq b$ de $X$. On a alors une suite exacte courte
			\begin{equation}\label{eq: une extension}
			0\To \H^1(X)\To \H^1(X\setminus \{a,b\}) \To \QQ(-1)\To 0
			\end{equation}
			qui exprime $\H^1(X\setminus \{a,b\})$ comme extension des groupes de cohomologie \enquote{purs} $\H^1(X)$ et $\QQ(-1)$, de poids respectifs $1$ et $2$. 
			
			Un point important est que l'extension \eqref{eq: une extension} est rarement scind\'{e}e, c'est-\`{a}-dire que $\H^1(X\setminus \{a,b\})$ n'est pas la somme directe $\H^1(X)\oplus \QQ(-1)$ -- en tout cas, pas de mani\`{e}re qui soit compatible aux isomorphismes de comparaison entre diff\'{e}rentes th\'{e}ories cohomologiques, ni \`{a} la th\'{e}orie de Hodge ou \`{a} l'action du groupe de Galois absolu.  Au niveau des p\'{e}riodes, l'obstruction \`{a} un tel scindage est mesur\'{e}e par des int\'{e}grales de formes diff\'{e}rentielles \enquote{de troisi\`{e}me esp\`{e}ce}, c'est-\`{a}-dire avec des p\^{o}les d'ordre $1$ en $a$ et $b$. Il est clair, cependant, que l'extension \eqref{eq: une extension} est scind\'{e}e si le diviseur $(a)-(b)$ de degr\'{e} $0$ est trivial dans le groupe de Picard $\operatorname{Pic}^0(X)$, c'est-\`{a}-dire s'il existe une fonction rationnelle sur $X$ avec seulement un z\'{e}ro d'ordre $1$ en $a$ et un p\^{o}le d'ordre $1$ en $b$. Cela sugg\`{e}re une relation entre la structure des extensions telles que \eqref{eq: une extension} et les cycles alg\'{e}briques.
				
			Cette id\'{e}e est d\'{e}velopp\'{e}e par Beilinson (influenc\'{e} par des id\'{e}es de Deligne), qui dans un article important \cite{beilinsonheight} postule en 1984 l'existence d'une cat\'{e}gorie ab\'{e}lienne de motifs (mixtes) $\mathsf{Mot}(k)$ qui contiendrait la cat\'{e}gorie $\mathsf{NumMot}(k)$ et o\`{u} les groupes d'extensions seraient reli\'{e}s aux cycles alg\'{e}briques. Par exemple, dans le cas d'une courbe alg\'{e}brique projective lisse $X$ on aurait :
			$$\operatorname{Ext}^1_{\mathsf{Mot}(k)}(\QQ(-1),\operatorname{M}^1(X)) \simeq \operatorname{Pic}^0(X)\otimes_\ZZ\QQ$$
			et l'extension \eqref{eq: une extension} correspondrait \`{a} la diff\'{e}rence $(a)-(b)\in\operatorname{Pic}^0(X)$.
			
			Beilinson va plus loin et \'{e}nonce une conjecture qui pr\'{e}voit que les valeurs sp\'{e}ciales des fonctions $L$ de motifs purs (et donc de toutes les fonctions $L$ classiques issues de l'arithm\'{e}tique) sont contr\^{o}l\'{e}es de mani\`{e}re pr\'{e}cise par les extensions entre motifs purs. La conjecture de Birch et Swinnerton-Dyer en est un cas particulier. Conclusion : m\^{e}me si l'on ne s'int\'{e}resse qu'aux vari\'{e}t\'{e}s alg\'{e}briques projectives et lisses, on ne peut se passer des motifs mixtes !
		
	\subsection{Motifs mixtes \`{a} la Voevodsky}\label{par: motifs mixtes}
	
		Il semble difficile d'\'{e}tendre la d\'{e}finition des motifs purs \`{a} la Grothendieck, qui ne concernent que les vari\'{e}t\'{e}s alg\'{e}briques projectives lisses, aux motifs mixtes, qui sont cens\'{e}s \^{e}tre une th\'{e}orie cohomologique universelle pour toutes les vari\'{e}t\'{e}s alg\'{e}briques. Un premier point de blocage est que \eqref{eq: morphisme corr} n'est pas d\'{e}fini si $Y$ n'est pas projective lisse. 
		
		Dans la vision de Grothendieck, la cat\'{e}gorie des motifs est \emph{ab\'{e}lienne}, avec des propri\'{e}t\'{e}s formelles proches de la cat\'{e}gorie des espaces vectoriels o\`{u} vivent les groupes de cohomologie. Beaucoup de progr\`{e}s de la th\'{e}orie des motifs \`{a} partir des ann\'{e}es 1980 vont reposer sur la suggestion, due \`{a} Deligne, de baser la th\'{e}orie des motifs non pas sur les groupes de cohomologie, mais sur les \emph{complexes} qui les calculent. Plut\^{o}t que de chercher \`{a} d\'{e}finir directement une cat\'{e}gorie ab\'{e}lienne des motifs, on cherche donc \`{a} 
		\begin{itemize}
		\item d'abord d\'{e}finir une cat\'{e}gorie \emph{triangul\'{e}e} des motifs, qui aurait les propri\'{e}t\'{e}s attendues de la cat\'{e}gorie d\'{e}riv\'{e}e de la cat\'{e}gorie ab\'{e}lienne des motifs\footnote{La cat\'{e}gorie d\'{e}riv\'{e}e $\operatorname{D}(\mathsf{A})$ d'une cat\'{e}gorie ab\'{e}lienne $\mathsf{A}$ s'obtient \`{a} partir de la cat\'{e}gorie des complexes dans $\mathsf{A}$ en inversant formellement les quasi-isomorphismes. Elle v\'{e}rifie les axiomes d'une cat\'{e}gorie triangul\'{e}e, d\'{e}finis par Verdier d'apr\`{e}s Grothendieck.} ;
		\item puis en extraire la cat\'{e}gorie ab\'{e}lienne des motifs gr\^{a}ce \`{a} un formalisme (\enquote{t-structures}) mis en place par Beilinson--Bernstein--Deligne--Gabber \cite{bbdg}.
		\end{itemize}
		
		Il y a de bonnes raisons de privil\'{e}gier le cadre triangul\'{e} au cadre ab\'{e}lien. D'abord, les cat\'{e}gories d\'{e}riv\'{e}es sont le r\'{e}ceptacle naturel des foncteurs d\'{e}riv\'{e}s, qui sont le langage de la cohomologie. En outre, on sait depuis les travaux de Grothendieck et de ses collaborateurs sur la cohomologie \'{e}tale $\ell$-adique que dans le cadre des familles de vari\'{e}t\'{e}s alg\'{e}briques les propri\'{e}t\'{e}s de la cohomologie s'expriment plus ais\'{e}ment dans la cat\'{e}gorie d\'{e}riv\'{e}e d'une cat\'{e}gorie de faisceaux (formalisme des \enquote{six op\'{e}rations}).  Par ailleurs, les groupes d'extensions entre motifs, centraux dans les visions de Beilinson et de Deligne, s'incarnent simplement dans le cadre triangul\'{e} comme des groupes de morphismes, ce qui rend tr\`{e}s coh\'{e}rent le slogan (M6) \enquote{les morphismes entre motifs sont reli\'{e}s aux cycles alg\'{e}briques}.
		
		La construction d'une cat\'{e}gorie triangul\'{e}e de motifs $\mathsf{DMot}(k)$ a \'{e}t\'{e} effectu\'{e}e par Voevodsky \cite{voevodskytriangulated} et ind\'{e}pendamment par Hanamura \cite{hanamuramixed} et Levine \cite{levinemixed}. Les morphismes dans $\mathsf{DMot}(k)$ sont reli\'{e}s aux cycles alg\'{e}briques via les groupes de Chow sup\'{e}rieurs de Bloch \cite{blochhigher} :
		$$\operatorname{Hom}_{\mathsf{DMot}(k)}(\QQ(-n),\operatorname{M}(X)[i]) \simeq \operatorname{CH}^n(X,2n-i).$$
		Les motifs de Chow $\mathsf{CHMot}(k)$ forment une sous-cat\'{e}gorie pleine de $\mathsf{DMot}(k)$.  Contrairement \`{a} la cat\'{e}gorie des motifs purs, la cat\'{e}gorie $\mathsf{DMot}(k)$ n'est cependant pas bas\'{e}e sur les groupes de Chow et leurs quotients mais sur le formalisme des \emph{correspondances finies} (ou \enquote{morphismes multivalu\'{e}s}) qui permettent de s'affranchir de subtilit\'{e}s li\'{e}es \`{a} la th\'{e}orie de l'intersection et notamment des \enquote{lemmes de d\'{e}placement} des cycles alg\'{e}briques.
		
		Il existe plusieurs d\'{e}finitions de $\mathsf{DMot}(k)$, plus ou moins pratiques en fonction des points de vue. Un ingr\'{e}dient central dans toutes les constructions est la \enquote{$\mathbb{A}^1$-localisation} : on force, pour toute vari\'{e}t\'{e} alg\'{e}brique $X$,  la projection $X\times\mathbb{A}^1\to X$ \`{a} induire un isomorphisme au niveau des motifs. Un tel isomorphisme existe en effet dans toute th\'{e}orie de cohomologie -- par exemple, en cohomologie de Betti, cela r\'{e}sulte du fait que $\mathbb{A}^1(\CC)=\CC$ est contractile.
		
		
		Il semble donc que la premi\`{e}re moiti\'{e} du programme soit accomplie, et il reste \`{a} extraire une cat\'{e}gorie ab\'{e}lienne de la cat\'{e}gorie $\mathsf{DMot}(k)$. Malheureusement, ce probl\`{e}me d'existence d'une \enquote{t-structure motivique} est toujours ouvert, et est probablement difficile puisqu'il implique les conjectures standard de Grothendieck ! 
		
		Alors, la th\'{e}orie des motifs est-elle condamn\'{e}e \`{a} \^{e}tre un grand \'{e}difice inutile tant qu'on n'a pas d\'{e}montr\'{e} des conjectures difficiles sur les cycles alg\'{e}briques ? La r\'{e}ponse est un grand non : le formalisme triangul\'{e} des motifs a amen\'{e} des r\'{e}sultats nouveaux, comme la conjecture de Bloch--Kato d\'{e}montr\'{e}e par Voevodsky \cite{voevodskyblochkato} en utilisant sa cat\'{e}gorie $\mathsf{DMot}(k)$. La fin de ce texte est consacr\'{e}e aux applications des motifs \`{a} l'\'{e}tude des p\'{e}riodes, qui reposent sur les aspects galoisiens des motifs.
		
\section{Aspects galoisiens et applications aux p\'{e}riodes}\label{par: aspects galoisiens}
		
	\subsection{Groupes de Galois motiviques}\label{par: tannaka}
	
		On cherche \`{a} construire une cat\'{e}gorie $\mathsf{Mot}(k)$ qui soit \emph{tannakienne} ; notamment elle doit \^{e}tre ab\'{e}lienne, mono\"{i}dale sym\'{e}trique rigide (existence de duaux), et tout foncteur de r\'{e}alisation $\mathsf{Mot}(k)\to \mathsf{Vect}_F$ est exact et fid\`{e}le (\emph{foncteur fibre}). Contentons-nous ici du cas \emph{neutre}, c'est-\`{a}-dire tel qu'il existe un foncteur fibre \`{a} valeurs dans les espaces vectoriels sur $F=\QQ$ (par exemple le foncteur de r\'{e}alisation de Betti pour $k$ un sous-corps de $\CC$). 
		
		Le formalisme tannakien, d\'{e}velopp\'{e} par Saavedra \cite{saavedra} et Deligne \cite{delignetannaka} d'apr\`{e}s des id\'{e}es de Grothendieck, est une \enquote{machine qui produit des groupes de Galois}. Pour un foncteur de r\'{e}alisation $\mathsf{Mot}(k)\to \mathsf{Vect}_\QQ$, on a un sch\'{e}ma en groupes\footnote{Un sch\'{e}ma en groupes d\'{e}fini sur $\QQ$ est la limite projective de groupes alg\'{e}briques d\'{e}finis sur $\QQ$, c'est-\`{a}-dire de sous-groupes de $\operatorname{GL}_N$ d\'{e}finis par des conditions polynomiales \`{a} coefficients dans $\QQ$ en les $N^2$ coefficients des matrices.} $G$ d\'{e}fini sur $\QQ$, appel\'{e} \emph{groupe de Galois motivique}, qui agit lin\'{e}airement sur la r\'{e}alisation de tout objet de $\mathsf{Mot}(k)$, et donc sur tous les groupes de cohomologie $\H^n(X)$ et plus g\'{e}n\'{e}ralement $\H^n(X,Y)$, de mani\`{e}re compatible \`{a} tous les morphismes dans $\mathsf{Mot}(k)$. Cette action est si riche qu'elle conna\^{i}t toute la structure des motifs : le foncteur naturel
		 $$\mathsf{Mot}(k) \longrightarrow \mathsf{Rep}(G)$$
		 induit par le foncteur de r\'{e}alisation est une \'{e}quivalence de cat\'{e}gories. On peut donc \'{e}tudier la structure des motifs via la th\'{e}orie des repr\'{e}sentations des groupes alg\'{e}briques. 
		 
		 Notons au passage que le formalisme tannakien sugg\`{e}re de \emph{d\'{e}finir} la cat\'{e}gorie des motifs comme la cat\'{e}gorie des repr\'{e}sentations d'un certain sch\'{e}ma en groupes (ou dualement comme la cat\'{e}gorie des comodules sur une alg\`{e}bre de Hopf). C'est l'approche suivie par Nori \cite{nori, hubermuellerstach} et par Ayoub \cite{ayoubgalois}
		 qui d\'{e}finissent deux candidats pour le groupe de Galois motivique. Bien que leurs approches soient fondamentalement diff\'{e}rentes (celle de Nori est relativement \'{e}l\'{e}mentaire alors que celle d'Ayoub repose sur le formalisme triangul\'{e} des motifs \`{a} la Voevodsky), on sait maintenant gr\^{a}ce aux travaux de Choudhury--Gallauer \cite{choudhurygallauer} qu'elles produisent les m\^{e}mes groupes de Galois motiviques. Dans le cas des motifs purs, ces groupes de Galois motiviques co\"{i}ncident avec ceux produits par la th\'{e}orie des cycles motiv\'{e}s d'Andr\'{e} \cite{andreinconditionnelle}.

		\subsection{Vers une th\'{e}orie de Galois des p\'{e}riodes}\label{par: galois periodes}
		
			Int\'{e}ressons-nous au cas $k=\QQ$, et appliquons le formalisme tannakien \`{a} une cat\'{e}gorie tannakienne de motifs $\mathsf{Mot}(k)$ et au foncteur de r\'{e}alisation de Betti. Le groupe de Galois motivique $G$ agit lin\'{e}airement sur tous les groupes de cohomologie $\H^n_{\B}(X)$ et plus g\'{e}n\'{e}ralement $\H^n_{\B}(X,Y)$. En se rappelant que les p\'{e}riodes apparaissent comme les coefficients de l'isomorphisme de comparaison \eqref{eq: B dR comp} entre Betti et de Rham, cela sugg\`{e}re que $G$ agit sur l'alg\`{e}bre des p\'{e}riodes via la formule :
			$$g \cdot \int_\sigma\omega = \int_{g\cdot\sigma}\omega.$$
			Probl\`{e}me : il se pourrait \emph{a priori} que deux \'{e}critures d'une m\^{e}me p\'{e}riode comme coefficient matriciel de \eqref{eq: B dR comp} donnent deux actions diff\'{e}rentes de $G$. La conjecture des p\'{e}riodes de Grothendieck pr\'{e}voit que ce ph\'{e}nom\`{e}ne ne se produit jamais. Pour l'\'{e}noncer, notons $T$ le torseur, produit par la \enquote{machine tannakienne}, des isomorphismes entre les foncteurs de r\'{e}alisation de de Rham et de Betti, et soit $\operatorname{comp}\in T(\CC)$ son point complexe induit par \eqref{eq: B dR comp}.

			\begin{conj}[Conjecture des p\'{e}riodes]
			Le point complexe $\operatorname{comp}$ est un point g\'{e}n\'{e}rique de $T$.
			\end{conj}
			
			Plus simplement, cette conjecture affirme que les p\'{e}riodes ne v\'{e}rifient pas de relations alg\'{e}briques qui ne soient pas \enquote{explicables} par des morphismes entre motifs.  Elle se sp\'{e}cialise en des conjectures de transcendance tr\`{e}s concr\`{e}tes sur certaines familles de p\'{e}riodes, rarement d\'{e}montr\'{e}es \cite{andregalois} mais qui peuvent \^{e}tre confirm\'{e}es exp\'{e}rimentalement. Notons tout de m\^{e}me les r\'{e}sultats g\'{e}n\'{e}raux de W\"{u}stholz \cite{wuestholzICM} et Huber--W\"{u}stholz \cite{huberwuestholz} au sujet des relations entre p\'{e}riodes de groupes de cohomologie de degr\'{e} $1$.
			
			La conjecture des p\'{e}riodes pour les motifs de Nori est \'{e}quivalente \cite{hubermuellerstach} \`{a} l'\'{e}nonc\'{e} \'{e}l\'{e}mentaire suivant \cite{kontsevichzagier}.
		
			\begin{conj}[Conjecture des p\'{e}riodes de Kontsevich--Zagier]
			Toutes les relations $\QQ$-lin\'{e}aires entre p\'{e}riodes sont des cons\'{e}quences de trois familles de relations :
			\begin{itemize}
			\item la bilin\'{e}arit\'{e} de l'int\'{e}gration ;
			\item la formule de changement de variable ;
			\item la formule de Stokes.
			\end{itemize}
			\end{conj}
			
			Si l'on croit \`{a} la conjecture des p\'{e}riodes on obtient, en suivant Andr\'{e} \cite{andregalois}, une th\'{e}orie de Galois des p\'{e}riodes, qui contient la th\'{e}orie de Galois classique des nombres alg\'{e}briques, et o\`{u} les groupes de Galois sont en g\'{e}n\'{e}ral des groupes alg\'{e}briques d\'{e}finis sur $\QQ$. Mais m\^{e}me sans savoir d\'{e}montrer la conjecture des p\'{e}riodes, le formalisme tannakien peut \^{e}tre utilis\'{e} pour d\'{e}montrer des r\'{e}sultats impressionnants sur les p\'{e}riodes, dont on donne maintenant deux exemples.

	\subsection{Applications des motifs aux p\'{e}riodes}\label{par: applications}
	
		\subsubsection{Valeurs z\^{e}ta multiples, et le th\'{e}or\`{e}me de Brown}

			Les \emph{valeurs z\^{e}ta multiples} sont des sommes de s\'{e}ries multiples qui g\'{e}n\'{e}ralisent les valeurs de la fonction z\^{e}ta de Riemann aux entiers :
			$$\zeta(n_1,\ldots,n_r)=\sum_{0 <  k_1< \cdots < k_r} \frac{1}{k_1^{n_1}\cdots k_r^{n_r}} ,$$
			o\`{u} les param\`{e}tres $n_i$ sont des entiers $\geq 1$ avec $n_r\geq 2$. Elles apparaissent naturellement dans plusieurs domaines des math\'{e}matiques et de la physique, notamment via le calcul d'int\'{e}grales de Feynman. Ce sont des p\'{e}riodes d'objets d'une cat\'{e}gorie tannakienne de motifs, les motifs de Tate mixtes sur $\ZZ$ -- un des cadres pour lesquels on sait produire une cat\'{e}gorie tannakienne \`{a} partir du formalisme triangul\'{e} de Voevodsky --, qui devrait former une sous-cat\'{e}gorie tannakienne de $\mathsf{Mot}(\QQ)$. Le qualificatif \enquote{de Tate mixte} d\'{e}signe les extensions it\'{e}r\'{e}es des motifs $\QQ(-n)=\QQ(-1)^{\otimes n}$ pour $n\in\ZZ$.
			
			Le th\'{e}or\`{e}me de Brown \cite{brownMTZ} est un r\'{e}sultat de structure sur la cat\'{e}gorie des motifs de Tate mixtes sur $\ZZ$, qui a comme cons\'{e}quence concr\`{e}te l'\'{e}nonc\'{e} suivant, conjectur\'{e} par Hoffman \cite{hoffman}.
			
			\begin{thm}
			Toute valeur z\^{e}ta multiple peut s'\'{e}crire comme combinaison $\QQ$-lin\'{e}aire des valeurs z\^{e}ta multiples $\zeta(n_1,\ldots,n_r)$ avec $n_i\in\{2,3\}$ pour tout $i$.
			\end{thm}
			
			L'ingr\'{e}dient central de la preuve est la th\'{e}orie de Galois (motivique) des valeurs z\^{e}ta multiples, initi\'{e}e par Goncharov \cite{goncharovgaloissymmetries}. Aussi surprenant que cela puisse para\^{i}tre, on ne conna\^{i}t pas \`{a} ce jour de preuve \enquote{non motivique} du th\'{e}or\`{e}me pr\'{e}c\'{e}dent. On renvoit le lecteur au livre de Burgos Gil et Fres\'{a}n \cite{burgosfresan} pour plus de d\'{e}tails sur sa preuve.
			
			Le th\'{e}or\`{e}me de Brown a des cons\'{e}quences hors des motifs et des valeurs z\^{e}ta multiples, gr\^{a}ce aux liens entre les motifs de Tate mixtes sur $\ZZ$ et les groupes de tresses ou encore les espaces de modules de courbes. \`{A} titre d'exemple, il est un des ingr\'{e}dients de la preuve d'un th\'{e}or\`{e}me r\'{e}cent de Chan--Galatius--Payne \cite{changalatiuspayne} sur la taille de la cohomologie de l'espace de modules $\mathcal{M}_g$ qui param\`{e}tre les courbes projectives et lisses de genre $g$.
			
		\subsubsection{S\'{e}ries de p\'{e}riodes, et le th\'{e}or\`{e}me d'Ayoub}
		
			Le th\'{e}or\`{e}me d'Ayoub \cite{ayoubKZ}, dont nous ne donnerons pas l'\'{e}nonc\'{e} ici, est un analogue \enquote{fonctionnel} de la conjecture des p\'{e}riodes de Kontsevich--Zagier o\`{u} les p\'{e}riodes sont remplac\'{e}es par des \enquote{s\'{e}ries de p\'{e}riodes}, qui sont certaines s\'{e}ries de Laurent en une variable
			$$\sum_{i\gg -\infty} \left(\int_{[0,1]^n}f_i\right) t^i.$$
			La preuve repose sur l'analyse de la structure des motifs sur $k=\CC(t)$ et du groupe de Galois motivique (de Nori ou d'Ayoub) correspondant.
			
			Ce qui est frappant dans le th\'{e}or\`{e}me d'Ayoub est qu'il s'agit d'un \'{e}nonc\'{e} g\'{e}n\'{e}ral de transcendance qui s'applique \`{a} toutes les p\'{e}riodes au sens fonctionnel, en une variable. Son pendant arithm\'{e}tique (les conjectures des p\'{e}riodes du \S\ref{par: galois periodes}) semble beaucoup plus difficile. Par exemple, la conjecture des p\'{e}riodes implique que $\zeta(3)$ n'est pas un multiple rationnel de $\pi^3$, ce qu'on ne sait pas d\'{e}montrer. La th\'{e}orie des motifs reste cependant un guide pr\'{e}cieux pour structurer les relations (ainsi que l'absence de relation) entre p\'{e}riodes.
			
		
		\subsubsection{Pour aller plus loin}
	
		Le lecteur int\'{e}ress\'{e} trouvera plus d'informations dans le livre de r\'{e}f\'{e}rence d'Andr\'{e} \cite{andrebook}. Le livre de Kahn \cite{kahnbook} est une introduction aux id\'{e}es motiviques bas\'{e}e sur les fonctions z\^{e}ta et $L$. L'article introductif de Kahn \cite{kahnlecondemaths} et le livre de Murre, Nagel et Peters \cite{murrenagelpeters} contiennent des d\'{e}tails sur les cat\'{e}gories de motifs purs. Le livre de Cisinski et D\'{e}glise \cite{cisinskideglisebook} \'{e}tudie les cat\'{e}gories triangul\'{e}es de motifs mixtes. Pour l'approche de Nori aux motifs et les liens avec les p\'{e}riodes, on pourra se r\'{e}f\'{e}rer au livre de Huber et M\"{u}ller-Stach \cite{hubermuellerstach}. Le volume collectif \cite{otherxups}, avec des contributions de Fres\'{a}n, Rivoal, et l'auteur, constitue une introduction \`{a} l'\'{e}tude de l'arithm\'{e}tique des p\'{e}riodes via les motifs.

\vspace{.3cm}		
		
		\subsubsection*{Remerciements}

Merci \`{a} Giuseppe Ancona, Yves Andr\'{e}, Damien Calaque, Ricardo Campos, Javier Fres\'{a}n, Charlotte Hardouin, Fran\c{c}ois L\^{e}, et un relecteur anonyme pour leurs commentaires sur une version pr\'{e}liminaire de cet article. Merci \`{a} Anthony Genevois pour avoir r\'{e}alis\'{e} les figures de l'article.

\vspace{.3cm}

\bibliographystyle{alpha-fr}
\bibliography{biblio}

\end{document}